\documentclass[mathematics,article,accept,moreauthors]{mdpi}

\usepackage{amssymb}
\usepackage{color}
\setitemize{parsep=6pt,itemsep=0pt,leftmargin=*,labelsep=5.5mm}
\setenumerate{parsep=6pt,itemsep=0pt,leftmargin=*,labelsep=5.5mm}
\setlist[description]{itemsep=0mm}

\newcommand{\seqnum}[1]{\href{http://oeis.org/#1}{\textcolor{black}{\underline{#1}}}}

\newcommand{\Pc}{{\mathbb P}_c}
\newcommand{\cH}{{\cal H}}
\newcommand{\bfprop}[1]{\smallskip\noindent{\bf #1}}

\newcommand{\mean}{\mathop{\mathrm{mean}}}
\newcommand{\lcm}{\mathop{\mathrm{LCM}}}

\newcommand{\li}{\mathop{\mathrm{li}}}
\newcommand{\Li}{\mathop{\mathrm{Li}}}
\newcommand{\Exp}{\mathop{\mathrm{Exp}}}
\newcommand{\Gumbel}{\mathop{\mathrm{Gumbel}}}

\firstpage{1}
\pubvolume{7}
\issuenum{5}
\articlenumber{400}
\pubyear{2019}
\copyrightyear{2019}
\history{This is arXiv:1901.03785. \ See the published article at \ https://doi.org/10.3390/math7050400 }

\Title{Predicting Maximal Gaps in Sets of Primes}

\Author{Alexei Kourbatov$^{1}$ and Marek Wolf\,$^{2}$}

\AuthorNames{Alexei Kourbatov and Marek Wolf}

\address{%
$^{1}$ \quad JavaScripter.net, 15127 NE 24th St., \#578, Redmond, WA 98052, USA; \ akourbatov@gmail.com \\
$^{2}$ \quad Faculty of Mathematics and Natural Sciences, Cardinal Stefan Wyszynski University, W\'oycickiego~1/3,~Bldg.\,21, PL-01-938 Warsaw, Poland; \ m.wolf@uksw.edu.pl }

\abstract{
Let $q>r\ge1$ be coprime integers. Let ${\mathbb P}_c={\mathbb P}_c(q,r,{\cal H})$ be an increasing sequence of primes $p$ satisfying two conditions:
(i) $p\equiv r$ (mod $q$)  \,{\it and}\,
(ii) $p$ starts a prime $k$-tuple with a given pattern ${\cal H}$.
Let $\pi_c(x)$ be the number of primes in ${\mathbb P}_c$ not exceeding $x$.
We heuristically derive formulas predicting the growth trend of the maximal gap $G_c(x)=\max_{p'\le x}(p'-p)$ between successive primes $p,p'\in{\mathbb P}_c$.
Extensive computations for primes up to $10^{14}$ show that a simple trend formula
$$G_c(x) \sim \frac{x}{\pi_c(x)}\cdot(\log \pi_c(x) + O_k(1))$$
works well for maximal gaps between initial primes of $k$-tuples with $k\ge2$
(e.g., twin primes, prime triplets, etc.) in residue class $r$ (mod $q$).
For $k=1$, however, a more sophisticated formula
$$G_c(x) \sim \frac{x}{\pi_c(x)}\cdot\big(\log{\frac{\pi_{c}^{2}(x)}{x}}+O(\log q)\big)$$
gives a better prediction of maximal gap sizes.
The latter includes the important special case of maximal gaps in the sequence of all primes ($k=1$, $q=2$, $r=1$).
The distribution of appropriately rescaled maximal gaps $G_c(x)$ is close to the Gumbel extreme value distribution.
Computations suggest that almost all maximal gaps satisfy a generalized strong form of Cram\'er's conjecture.
We~also conjecture that the number of maximal gaps between primes in ${\mathbb P}_c$ below $x$ is $O_k(\log x)$.
}

\begin{document}


\section{Introduction}

A {\em prime gap} is the difference between consecutive prime numbers.
The sequence of prime gaps behaves quite erratically (see OEIS {\seqnum{A001223} \cite{oeis}).  
While the prime number theorem tells us that the {\em average} gap between primes near $x$ is about $\log x$,
the actual gaps near $x$ can be significantly larger or smaller than $\log x$.
We call a gap {\it maximal} if it is strictly greater than all gaps before it.
Large gaps between primes have been studied by many authors; see, e.g.,
\cite{ramanujannotebooks4,bhp2001,cramer,fgkmt,funkhouser2018,granville,nicely,nicelynyman2003,toes2014,shanks}.
In the early 1910s, Ramanujan considered maximal prime gaps up to low 7-digit primes \cite[p.\,133]{ramanujannotebooks4}.
More than a century later, we~know all  maximal gaps between primes below $2^{64}$  \cite{nicely}.

Let $G(x)$ be the maximal gap between primes not exceeding $x$:
$$
G(x) ~=~ \max_{p_{n+1}\le x} (p_{n+1} - p_n).
$$
Estimating $G(x)$ is a subtle and delicate problem.
Cram\'er \cite{cramer} conjectured on probabilistic grounds that $G(x)=O(\log^2 x)$,
while Shanks \cite{shanks} heuristically found that $G(x)\sim \log^2 x$.
Granville~\cite{granville} heuristically argued that for a certain subsequence of maximal gaps
we should expect significantly larger sizes of $G(x)$; namely,
$\limsup\limits_{x\to\infty}\frac{G(x)}{\log^2 x}\ge2e^{-\gamma}\approx1.1229$.

Baker, Harman, and Pintz \cite{bhp2001} proved that $G(x)=O(x^{0.525})$; indeed,
computation suggests that $G(x)<x^{0.525}$ for $x\ge153$. 
Ford, Green, Konyagin, Maynard, and Tao \cite{fgkmt} proved that the order of
$G(x)$ is at least $\frac{c \log x \log\log x\log\log\log\log x}{\log\log\log x}$, solving a long-standing conjecture of Erd\H{o}s.

\smallskip
Earlier, we independently proposed formulas closely related to the Cram\'er and Shanks conjectures.
Wolf \cite{wolf1998, wolf2011, wolf2014} expressed the probable size of maximal gaps $G(x)$
in terms of the prime-counting function~$\pi(x)$:
\begin{equation}\label{wolfGx}
G(x) ~\sim~ \frac{x}{\pi(x)} \cdot \left( \log\frac{\pi^2(x)}{x} + O(1) \right),
\end{equation}
which suggests an analog of Shanks conjecture $G(x)\sim \log^2 x - 2\log x \log\log x + O(\log x)$;
see also Cadwell \cite{cadwell}.
Extending the problem statement to {\em prime $k$-tuples}, Kourbatov \cite{kourbatov2013,kourbatov2013tables}
empirically tested (for $x\le10^{15}$, $k\le7$)
the following heuristic formula for the probable size of maximal gaps $G_k(x)$
between prime $k$-tuples below $x$:
\begin{equation}\label{kourbGkx}
G_k(x) ~\sim~ a(x)\cdot \left( \log\frac{x}{a(x)} + O(1)  \right),
\end{equation}
where $a(x)$ is the expected average gap between the particular prime $k$-tuples near $x$.
Similar to (\ref{wolfGx}), formula (\ref{kourbGkx}) also suggests an analog of the Shanks conjecture,
$G_k(x) \sim C \log^{k+1}x$, with a negative correction term of size
$O_k\big((\log x)^k \log\log x\big)$; see also \cite{ford2018,toes2015}.

In this paper we study a further generalization of the prime gap growth problem, viz.:
What~happens to maximal gaps if we only look at primes in a specific {\em residue class} mod $q$?
The~new problem statement subsumes, as special cases,
{\it maximal prime gaps} ($k=1$, $q=2$) as well as maximal gaps between {\it prime $k$-tuples} ($k\ge2$, $q=2$).
One of our present goals is to generalize formulas (\ref{wolfGx}) and~(\ref{kourbGkx}) to gaps between primes
in a residue class---and test them in computational experiments.
Another goal is to investigate {\em how many} maximal gaps should be expected between primes $p\le x$ in a residue class,
with an additional (optional) condition that $p$ starts a prime constellation of a certain type.

\subsection{Notation}
%
\noindent
\begin{tabular}{@{}ll}
$q$, $r$          & coprime integers, $1\le r<q$                                                           \\
$p_n$             & the $n$-th prime;\, $\{p_n\} = \{2,3,5,7,11,\ldots\}$                                  \\
$\Pc=\Pc(q,r,\cH)$ & increasing sequence of primes $p$ such that (i) $p\equiv r$ (mod $q$) and              \\
& (ii) $p$ is the least prime in a prime $k$-tuple with a given pattern $\cH$.        \\
& {\em Note:} ${\mathbb P}_c$ depends on $q$, $r$, $k$, and on the pattern $\cH$ of the $k$-tuple.    \\
& When $k=1$, \ ${\mathbb P}_c$ is the sequence of {\em all} primes $p\equiv r$ (mod $q$).           \\
$\cH$  & the $k$-tuple pattern of offsets: $\cH = (\Delta_1,\Delta_2,\ldots,\Delta_k)$ (see Section \ref{definitions}) \\
$\gcd(m,n)$       & the greatest common divisor of $m$ and $n$                                             \\
$\varphi(q)$      & Euler's totient function (OEIS \seqnum{A000010})                                       \\
$\varphi_{k,\cH}(q)$ & Golubev's generalization (\ref{eqdefphik}) of Euler's totient (see Section \ref{defphik})         \\
$\Gumbel(x;\alpha,\mu)$& the Gumbel distribution cdf: \ $\Gumbel(x;\alpha,\mu) =
e^{-e^{-{x-\mu\over\vphantom{f}\alpha}}}$                                              \\
$\Exp(x;\alpha)$  & the exponential distribution cdf: \ $\Exp(x;\alpha)=1-e^{-x/\alpha}$                   \\
$\alpha$          & the {\em scale parameter} of exponential/Gumbel distributions, as applicable           \\
$\mu$             & the {\em location parameter} ({\em mode}) of the Gumbel distribution                   \\
$\gamma$          & the Euler--Mascheroni constant: \ $\gamma = 0.57721\ldots$                             \\
$C_k=C_{k,\cH}$    & the Hardy--Littlewood constants (see {Appendix} \ref{hlconst})                     \\
$\log x$          & the natural logarithm of $x$                                                           \\
$\li x$           & the logarithmic integral of $x$: \
$\displaystyle\li x \,= \int_0^x{\negthinspace}{dt\over\log t}
\,= \int_2^x{\negthinspace}{dt\over\log t} + 1.04516\ldots$        \\
$\Li_k(x)$        & the integral $\displaystyle\int_2^x{\negthinspace}{dt\over\log^k t}$
(see {Appendix} \ref{integralsLikx})           \\
\end{tabular}

\noindent
\begin{tabular}{@{}ll}
$\phantom{1234567890123}$
& {\em Gap measure functions:}                             ${\large\phantom{1^{1^1}}}$   \\
$G(x)$            & the maximal gap between primes $\le x$                                                 \\
$G_{q,r}(x)$      & the maximal gap between primes $p=r+nq \le x$                (case $k=1$)              \\
$G_{c}(x)$        & the maximal gap between primes $p\in{\mathbb P}_c$ not exceeding $x$                   \\
$R_{c}(n)$        & the $n$-th record (maximal) gap between primes $p\in{\mathbb P}_c$                     \\
$a$, $ a_c$, $\bar a_c$
& the expected average gaps between primes in ${\mathbb P}_c$ (see Section~\ref{avgaps})  \\
$T$, $ T_c$, $\bar T_c$
& trend functions predicting the growth of maximal gaps (see Section~\ref{maxgaps})
\\
& {\em Gap counting functions:}                            ${\large\phantom{1^{1^1}}}$   \\
$N_{c}(x)$        & the number of maximal gaps $G_{c}$ with endpoints $p\le x$                             \\
$N_{q,r}(x)$      & the number of maximal gaps $G_{q,r}$ with endpoints $p\le x$ (case $k=1$)              \\
$\tau_{q,r}(d,x)$ & the number of gaps of a given even size $d=p'-p$ between successive                    \\
& primes $p,p'\equiv r$ (mod $q$), with $p'\le x$; \ \
$\tau_{q,r}(d,x)=0$ \ if $q\nmid d$ or $2\nmid d$.
\\
& {\em Prime counting functions:}                          ${\large\phantom{1^{1^1}}}$   \\
$\pi(x)$          & the total number of primes $p_n\le x$                                                  \\
$\pi_{c}(x)$      & the total number of primes $p\in{\mathbb P}_c$ not exceeding $x$                       \\
$\pi(x;q,r)$      & the total number of primes $p=r+nq\le x$ (case $k=1$)                                  \\
\end{tabular}


\vspace{12pt}
Quantities with the $c$ subscript may, in general, depend on
$q$, $r$, $k$, and on the pattern of the prime $k$-tuple.
However, average gaps $a$, $ a_c$, $\bar a_c$ and
trend functions $T$, $ T_c$, $\bar T_c$ are independent of $r$.
Expressions like $\pi_c^2(x)$ or $\log^2 x$ denote the {\em square} of the respective function.

\subsection{Definitions: Prime $k$-Tuples, Gaps, Sequence ${\mathbb P}_c$}\label{definitions}

{\em Prime $k$-tuples} are clusters of $k$ consecutive primes that have an admissible\footnote{
A $k$-tuple is admissible (infinitely repeatable) unless it is prohibited by an elementary divisibility argument.
For~example, the cluster of five numbers ($p$,~$p+2$, $p+4$, $p+6$, $p+8$) is prohibited because
one of the numbers is divisible by 5 (and, moreover, {\it at least one} of the numbers is divisible by 3);
hence all these five numbers cannot simultaneously be prime infinitely often.
Likewise, the cluster of three numbers ($p$,~$p+2$, $p+4$) is prohibited because
one of the numbers is divisible by 3; so these three numbers cannot simultaneously be prime infinitely often.
} pattern $\cH$.
%
In what follows, when we speak of a $k$-tuple, for~certainty we will mean a {\it densest admissible prime $k$-tuple},
with a given $k\le7$. However, our observations can be extended to other admissible $k$-tuples,
including those with larger $k$ and not necessarily densest ones.
The densest $k$-tuples that exist for a given $k$ may sometimes be called
{\em prime constellations} or {\em prime $k$-tuplets}.
Below are examples of prime $k$-tuples with $k=2$, 4, 6.
\begin{itemize}
\item
{\em Twin primes} are pairs of consecutive primes that have the form ($p$, $p+2$).
This is the densest admissible pattern of two; 
$\cH=(0,2)$.

\item
{\em Prime quadruplets} are clusters of four consecutive primes of the form ($p$, $p+2$, $p+6$, $p+8$).
This~is the densest admissible pattern of four;
$\cH=(0,2,6,8)$.

\item
{\em Prime sextuplets} are clusters of six consecutive primes 
($p$, $p+4$, $p+6$, $p+10$, $p+12$, $p+16$).
This~is the densest admissible pattern of six;
$\cH=(0,4,6,10,12,16)$.
\end{itemize}

A {\em gap} between prime $k$-tuples is the
distance $p'-p$ between the initial primes $p$ and $p'$ in two consecutive $k$-tuples of the same type
(i.e., with the same pattern).
For example, the gap between twin prime pairs $(17,19)$ and $(29,31)$ is 12: \
$p'-p ~=~ 29 - 17 ~=~ 12.$

A {\em maximal gap} between prime $k$-tuples is a gap that is strictly greater than all gaps
between preceding $k$-tuples of the same type.
For example, the gap of size 6 between twin primes $(5,7)$ and $(11,13)$ is maximal, while
the gap (also of size 6) between twin primes $(11,13)$ and $(17,19)$ is {\bf\em not}~maximal.

\pagebreak
Let $q>r\ge1$ be coprime integers.
Let $\Pc=\Pc(q,r,\cH)$ be an increasing sequence of primes $p$ satisfying two conditions:
(i) $p\equiv r$ (mod $q$) and (ii) $p$ starts a prime $k$-tuple with a given pattern $\cH$.
Importantly, $\Pc$ depends on $q$, $r$, $k$, and on the pattern of the $k$-tuple.
When $k=1$, \ $\Pc$ is the sequence of {\em all} primes $p\equiv r$ (mod $q$).
Gaps between primes in $\Pc$ are defined as differences $p'-p$
between successive primes $p,p'\in\Pc$.
As before, a gap is maximal if it is strictly greater than all preceding gaps.
Accordingly, for successive primes $p,p'\in\Pc$ we define
%
$$
G_c(x) = \max\limits_{p'\le\,x \atop p,p'\in\,\Pc} (p'-p).
$$

Studying maximal gaps between primes in ${\mathbb P}_c$ is convenient. Indeed, if the modulus $q$
used for defining ${\mathbb P}_c$ is ``not too small'', we get {\em plenty of data} to study maximal gaps;
that is, we get many sequences of maximal gaps corresponding to ${\mathbb P}_c$'s with different $r$ for the same $q$,
which allows us to study common properties of these sequences. (One such property is the
{\em average} number of maximal gaps between primes in ${\mathbb P}_c$ below $x$.)
By contrast, data on maximal prime gaps are scarce:
at present we know that there are only 80 maximal gaps between primes below $2^{64}$ \cite{nicely}.
Even fewer maximal gaps are known between $k$-tuples of any given type \cite{kourbatov2013tables}.

\begin{Remark}~~~~~~~~~~~~~~~~~~

\begin{enumerate}[align=parleft,leftmargin=*,labelsep=3mm]
\item[(i)] In Section \ref{heuristicsection} we derive formulas predicting the {\it most probable sizes} of maximal gaps $G_c(x)$. 
It is not known how close these most probable sizes might be to the {\it maximal order}
of $G_c(x)$.  
Thus, in the special case $k=1$, $q=2$, probable values of $G(x)$ seem to be $\sim\log^2 x - 2\log x\log\log x$ \cite{wolf2011};
but it is not implausible that the maximal order of $G(x)$ is closer to $2e^{-\gamma}\log^2 x$ \ \cite{granville}.
For further discussion of extremely large gaps, see Section~\ref{ExtraLargeGaps}.


\item[(ii)]How hard is it to compute gaps in sequence $\Pc$?
Given $k=1$, $q\approx10^3$ and $r$ coprime to $q$, our PARI/GP code (Appendix \ref{AppendixA})
takes several hours to compute all maximal gaps in sequence $\Pc$ up to 14-digit primes.
In~some numerical experiments, we carried out the computation all the way to $10^{14}$.
In most cases, however, we stopped the computation at $e^{28}$ or at $10^{12}$ or even earlier, to quickly gather statistics for all $r$ coprime to $q$.
A~similar strategy was also used for sequences $\Pc$ with $k\ge2$ (source code for $k\ge2$ is not included).
See~Section~\ref{numresults} for a detailed discussion of our numerical results.
\end{enumerate}
\end{Remark}

\subsection{Generalization to Other Subsets of Primes}\label{othersubsets}
Sequences $\Pc$ include, as special cases, many different subsets of prime numbers: primes in a given residue class, twin primes, triplets, quadruplets, etc.
However, formulas akin to (\ref{wolfGx}) and (\ref{kourbGkx}) definitely have an even wider area of applicability.
Namely, we expect that certain analogs of (\ref{wolfGx}) or (\ref{kourbGkx}), possessing the general form
$$
\mbox{ maximal gap size } ~\sim~
(\mbox{average gap near }x) ~\cdot~ L(x), \quad
\mbox{ with } L(x) \lesssim c\log x,
$$
will also be applicable to maximal gaps in the following subsets of primes:

\begin{itemize}
\item the sequence of prime-indexed primes \cite[]{bb2009},  \seqnum{A006450}

\item higher iterates of prime-indexed primes \cite[]{bko2013,batchko2014,guariglia2019}, \seqnum{A038580}

\item primes $p=n^2+1$, \ $n\in{\mathbb N}$ \cite[]{wolf2013}, \seqnum{A002496}

\item primes $p=f(n)$, where $f(n)$ is an irreducible polynomial in $n$,

\item primes in sequences of Beatty type: $p=\lfloor \beta n + \delta \rfloor$, \ $n\in{\mathbb N}$, \
for a fixed irrational $\beta>1$ and a fixed real $\delta$ \cite[]{bakerzhao2016}, \seqnum{A132222}.
\end{itemize}

The above list is by no means exhaustive, but it may serve as a starting point for future work.

\pagebreak

\subsection{When Are Equations (1), (2) Inapplicable?}\label{inapplicability}

Analogs of Equations (\ref{wolfGx}) and (\ref{kourbGkx}) are {\bf\em not} applicable to sequences where (almost) every gap is maximal.
Examples of this kind include:

\begin{itemize}
\item Mills primes \cite[]{mills}, \seqnum{A051254},

\item base-$B$ repunit primes \cite[]{salas2011}, \seqnum{A076481},

\item  primes nearest to $e^n$  (\seqnum{A037028}),

\item in general, any sequence whose terms grow exponentially or super-exponentially.

\end{itemize}

\section{Heuristics and Conjectures}\label{heuristicsection}
We now focus on deriving  analogs of formulas (\ref{wolfGx}) and (\ref{kourbGkx}) for sequences $\Pc=\Pc(q,r,\cH)$.

\subsection{Equidistribution of $k$-Tuples}\label{basicconj}
Everywhere we assume that $q>r$ are coprime positive integers.
Let $\pi(x;q,r)$ be the number of primes $p\equiv r$ (mod $q$) such that $p\le x$.
The prime number theorem for arithmetic progressions \cite{fine2016,fgoldston1996} establishes that
\begin{equation}\label{PNTprog}
\pi(x;q,r) ~\sim~ \frac{\li x}{\varphi(q)} \qquad\mbox{ as }x\to\infty.
\end{equation}

Furthermore, the generalized Riemann hypothesis (GRH) implies that
\begin{equation}\label{GRH}
\pi(x;q,r) ~=~ \frac{\li x}{\varphi(q)}+O_\varepsilon(x^{1/2+\varepsilon}) \qquad\mbox{ for any }\varepsilon>0.
\end{equation}

That is to say, the primes below $x$ are approximately equally distributed among the
$\varphi(q)$ ``allowed'' residue classes (these classes form the {\it reduced residue system} modulo $q$).
Roughly speaking, the~GRH implies that, as $x\to\infty$, the numbers $\pi(x;q,r)$ and
$\lfloor \li x / \varphi(q)\rfloor$ almost agree in the left half of their~digits.

Based on empirical evidence, below we conjecture that a similar phenomenon also occurs for prime $k$-tuples:
in every {\em $\cH$-allowed} residue class (as defined in Section~\ref{defphik}) there are infinitely many primes
starting an admissible $k$-tuple with a particular pattern $\cH$.
Moreover, such primes are distributed approximately equally among all $\cH$-allowed residue classes modulo $q$.
Our conjectures are closely related to the Hardy--Littlewood $k$-tuple conjecture \cite{hl1923}
and the Bateman--Horn conjecture \cite{bh1962}.

\subsubsection{Counting the $\cH$-Allowed Residue Classes}\label{defphik}

Consider an example: take $\cH=(0,2)$.
Which residue classes modulo 4 may contain the lesser prime $p$ in a pair of twin primes $(p,\,p+2)$?
Clearly, the residue class 0 mod 4 is prohibited: all numbers in this class are even.
The residue class 2 mod 4 is prohibited for the same reason.
The remaining residue classes, $p\equiv1$ mod 4 and $p\equiv3$ mod 4, are not prohibited.
We call these two classes {\em $\cH$-allowed}. Indeed, each of these two residue classes
does contain lesser twin primes---and there are, conjecturally, infinitely many such primes
in each class (see OEIS \seqnum{A071695} and \seqnum{A071698}).

In general, given an admissible $k$-tuple with pattern $\cH=(\Delta_1,\Delta_2,\ldots,\Delta_k)$,
we say that a residue class $r$ (mod $q$) is {\it $\cH$-allowed} if
$$
\gcd(r+\Delta_1,q)=\gcd(r+\Delta_2,q)=\gcd(r+\Delta_3,q)=\ldots=\gcd(r+\Delta_k,q)=1.
$$

Thus a residue class is $\cH$-allowed if it is not prohibited
(by divisibility considerations) from containing infinitely many
primes $p$ starting a prime $k$-tuple with pattern $\cH$.

How many residue classes modulo $q$ are $\cH$-allowed?
To count them, we will need an appropriate generalization of Euler's totient $\varphi(q)$:
Golubev's totient functions \cite{golubev1953,golubev1958,golubev1962};
see also \cite[p.\,289]{SandorCrstici2004II}.


\begin{Definition}
{\bf\em Golubev's totient} $\varphi_{k,\cH}(q)$ is the number of $\cH$-allowed residue classes modulo $q$
for a given pattern $\cH=(\Delta_1,\ldots,\Delta_k)$.
More formally,
\begin{equation}\label{eqdefphik}
\varphi_{k,\cH}(q) ~=~
\sum\limits_{1\le x\le q\atop \gcd(x+\Delta_1,q)\,=\cdots=\,\gcd(x+\Delta_k,q)\,=\,1} 1.
\end{equation}
\end{Definition}

\begin{Example}
For prime quadruplets $(p, p+2, p+6, p+8)$ we have
$$
k=4,
\quad \cH = (\Delta_1, \Delta_2, \Delta_3, \Delta_4) = (0,2,6,8),
\quad\mbox{ and }
\quad \varphi_{4,\cH}(q)=\mbox{\seqnum{A319516}}(q).
$$
\end{Example}

\noindent
For instance, when $q=30$, we have $\varphi_{4,\cH}(q)=1$: indeed, there is only one residue class,
namely, $p\equiv11$ (mod 30) where divisibility considerations allow infinitely many primes $p$
at the beginning of prime quadruplets $(p, p+2, p+6, p+8)$.

Note that $\varphi_1(q)=\varphi(q)$ is Euler's totient function, \seqnum{A000010};
and, for densest admissible $k$-tuples,
$\varphi_{2,\cH}(q)$ is \seqnum{A002472}, see also Alder \cite{alder1958};
$\varphi_{3,\cH}(q)$ is \seqnum{A319534};
$\varphi_{4,\cH}(q)$ is \seqnum{A319516};
$\varphi_{5,\cH}(q)$ is \seqnum{A321029}; and
$\varphi_{6,\cH}(q)$ is \seqnum{A321030}.
Like Euler's totient, the functions $\varphi_{k,\cH}$ are {\it multiplicative} \cite{golubev1958}.

\subsubsection{The $k$-Tuple Infinitude Conjecture}\label{infinitude}

We expect each of the $\cH$-allowed residue classes $r$ (mod $q$) to contain infinitely many primes $p$
starting admissible prime $k$-tuples
with pattern $\cH$.
In other words, the corresponding sequence $\Pc=\Pc(q,r,\cH)$ is {\em infinite}.

\begin{Remark} ~~~~~~~~~

\begin{enumerate}[align=parleft,leftmargin=*,labelsep=3mm]
\item[(i)] The $k$-tuple infinitude conjecture generalizes Dirichlet's theorem on arithmetic progressions \cite{dirichlet}.
\item[(ii)] The conjecture follows from the Bateman--Horn conjecture \cite{bh1962}.
\end{enumerate}
\end{Remark}

\subsubsection{The $k$-Tuple Equidistribution Conjecture}

Suppose the residue class $r$ (mod $q$) is $\cH$-allowed.
We conjecture that the number of primes $p\in\Pc(q,r,\cH)$, $p\le x$, is
\begin{equation}\label{ktuplequidist}
\pi_c(x) ~=~ \frac{C_{k,\cH}}{\varphi_{k,\cH}(q)} \,{\Li}_k(x) + O_{\eta,\cH}(x^\eta) \quad\mbox{ as }x\to\infty,
\end{equation}

\noindent
where $\eta<1$, the coefficient $C_{k,\cH}$ is the {\it Hardy--Littlewood constant}
for the particular $k$-tuple ({Appendix}~\ref{hlconst}),
${\Li}_k(x)=\int_2^x \log^{-k}t\,dt$
({Appendix} \ref{integralsLikx}), and
$\varphi_{k,\cH}(q)$ is Golubev's totient function (\ref{eqdefphik}).

\begin{Remark} ~~~~~~~~~~~~

\begin{enumerate}[align=parleft,leftmargin=*,labelsep=3mm]
\item[(i)] Conjecture (\ref{ktuplequidist}) is akin to the GRH-based Equation~(\ref{GRH}); the latter pertains to the case $k=1$.
\item[(ii)] The conjecture is compatible with the Bateman--Horn and Hardy--Littlewood $k$-tuple conjectures
but does not follow from them.
\item[(iii)] It is plausible that, similar to (\ref{GRH}), in (\ref{ktuplequidist}) we can take $\eta=\frac{1}{2}+\varepsilon$
for any $\varepsilon>0$.
\end{enumerate}
\end{Remark}


\subsection{Average Gap Sizes}\label{avgaps}

Consider a sequence $\Pc=\Pc(q,r,\cH)$, where the residue class $r$ (mod $q$) is $\cH$-allowed.
We define the {\em expected average gaps} between primes in $\Pc$ as follows.


\begin{Definition}
The {\bf\em expected average gap} between primes in $\Pc$ {\bf\em below} $x$ is
\begin{equation}\label{tildea}
a_c(x) ~=~ \frac{\varphi_{k,\cH}(q)}{C_{k,\cH}} \cdot \frac{x}{\Li_k(x)}.
\end{equation}
\end{Definition}


\begin{Definition}
The {\bf\em expected average gap} between primes in ${\mathbb P}_c$ {\bf\em near} $x$ is
\begin{equation}\label{bara}
\bar a_c(x) ~=~ \frac{\varphi_{k,\cH}(q)}{C_{k,\cH}} \cdot \log^k x.
\end{equation}
\end{Definition}

In view of the equidistribution conjecture (\ref{ktuplequidist}),
it is easy to see from these definitions that
$$
\frac{x}{\pi_c(x)} ~\approx~
a_c(x) ~<~ \bar a_c(x) \qquad\mbox{ for large } x.   
$$

We have the limits (with very slow convergence):
\begin{equation}
\lim_{x\to\infty} \frac{ a_c(x)}{\bar a_c(x)} ~=~ 1,
\end{equation}
\begin{equation}
\lim_{x\to\infty} \frac{\bar a_c(x) -  a_c(x) }{ \bar a_c(x)}\cdot \log x ~=~ k.
\end{equation}

\subsection{Maximal Gap Sizes}\label{maxgaps}

Recall that formula (\ref{wolfGx}) is applicable to the special case
$q=2, k=1$ \cite{wolf1998, wolf2011, wolf2014},
while (\ref{kourbGkx}) is applicable to the special cases
$q=2, k\ge2$ \cite{kourbatov2013}.
We are now ready to generalize (\ref{wolfGx}) and (\ref{kourbGkx})
for predicting maximal gaps between primes in sequences ${\mathbb P}_c$ with $q\ge2$.

\subsubsection{Case of $k$-Tuples: $k\ge2$}
Consider a probabilistic example.
Suppose that intervals between rare random events are exponentially distributed, with cdf
$\Exp(\xi;\alpha)=1-e^{-\xi/\alpha}$, where $\alpha$ is the mean interval between events.
If our observations of the events continue for $x$ seconds,
extreme value theory (EVT) predicts that the expected maximal interval between events is
\begin{equation}\label{EVTtrend}
\mbox{ expected maximal interval}
~=~ \alpha \log\frac{x}{\alpha}+O(\alpha)
~=~ \frac{x}{\Pi(x)} \log\Pi(x)+O(\alpha),
\end{equation}

\noindent
where $\Pi(x)\approx x/\alpha$ is the total count of the events we observed in $x$ seconds.
(For details on deriving Equation~(\ref{EVTtrend}),
see e.g.~\cite[pp.\,114--116]{gumbel1958}  or \cite[Sect.\,8]{kourbatov2013}.)

By analogy with EVT, we define the expected trend functions for maximal gaps as follows.


\begin{Definition}
The {\bf\em lower trend} of maximal gaps between primes in ${\mathbb P}_c$ is
\begin{equation}\label{lowertrend}
T_c(x) ~=~  a_c(x) \cdot \log\frac{C_{k,\cH} \Li_k(x) }{ \varphi_{k,\cH}(q)} .
\end{equation}
\end{Definition}

In view of the equidistribution conjecture (\ref{ktuplequidist}),
\begin{equation}
T_c(x) ~\approx~  a_c(x)\cdot\log\pi_c(x)
~\approx~ \frac{x}{\pi_c(x)}\cdot\log\pi_c(x)
\qquad\mbox{ as }x\to\infty.
\end{equation}

We also define another trend function, $\bar T_c(x)$, which is simpler because it does not use $\Li_k(x)$.


\begin{Definition}
The {\bf\em upper trend} of maximal gaps between primes in ${\mathbb P}_c$ is
\begin{equation}\label{uppertrend}
\bar T_c(x) ~=~ \bar a_c(x) \cdot \log\frac{x }{ \bar a_c(x)} .
\end{equation}
\end{Definition}

The above definitions imply that
\begin{equation}
T_c(x) ~<~ \bar T_c(x) ~<~ {C_{k,\cH}^{-1}\, \varphi_{k,\cH}(q)} \cdot \log^{k+1}x
\qquad\mbox{ for large }x.
\end{equation}

At the same time, we have the asymptotic equivalence:

\begin{equation}
T_c(x) ~\sim~ \bar T_c(x) ~\sim~ {C_{k,\cH}^{-1}\, \varphi_{k,\cH}(q)} \cdot \log^{k+1}x
\qquad\mbox{ as }x\to\infty.
\end{equation}

We have the limits (convergence is quite slow):

\begin{equation}\label{trenddiff}
\lim_{x\to\infty} \frac{\bar T_c(x) -  T_c(x) }{ \bar a_c(x)} = k,
\end{equation}
\begin{equation}
\lim_{x\to\infty} \frac{C_{k,\cH}^{-1}\, \varphi_{k,\cH}(q) \log^{k+1} x - \bar T_c(x) }{\bar a_c(x)\log\log x} = k.
\end{equation}

Therefore, $\bar T_c(x) -  T_c(x) = O_k(\bar a_c)$, while \
$C_{k,\cH}^{-1}\, \varphi_{k,\cH}(q) \log^{k+1} x - \bar T_c(x) = O_k(\bar a_c \log\log x)$.

We make the following conjectures regarding the behavior of maximal gaps $G_c(x)$.

\vspace{6pt}
\bfprop{Conjecture on the trend of $G_c(x)$.}
For any sequence ${\mathbb P}_c$ with $k\ge2$, a positive proportion of maximal gaps $G_c(x)$
satisfy the double inequality

\begin{equation}\label{GcTrend}
T_c(x) ~\lesssim~ G_c(x) ~\lesssim~ \bar T_c(x)
\qquad\mbox{ as }x\to\infty,
\end{equation}

\noindent
and the difference $G_c(x)-\bar T_c(x)$ changes its sign infinitely often. 

\vspace{6pt}
\bfprop{Generalized Cram\'er conjecture for $G_{c}(p)$.}
Almost all maximal gaps $G_{c}(p)$ satisfy

\begin{equation}\label{GcCramer}
G_{c}(p) ~<~ C_{k,\cH}^{-1}\, \varphi_{k,\cH}(q) \log^{k+1} p.
\end{equation}

\vspace{6pt}
\bfprop{Generalized Shanks conjecture for $G_{c}(p)$.}
Almost all maximal gaps $G_{c}(p)$ satisfy
\begin{equation}\label{GcShanks}
G_{c}(p) ~\sim~ C_{k,\cH}^{-1}\, \varphi_{k,\cH}(q) \log^{k+1} p \qquad\mbox{ as }p\to\infty.
\end{equation}
Here $G_{c}(p)$ denotes the maximal gap that ends at the prime $p$.

\subsubsection{Case of Primes: $k=1$}\label{caseOfPrimes}

The EVT-based trend formulas (\ref{lowertrend}) and (\ref{uppertrend})
work well for maximal gaps between $k$-tuples, $k\ge2$. However, when $k=1$,
the observed sizes of maximal gaps $G_{q,r}(x)$ between primes in residue class $r$ mod $q$
are usually a little less than predicted by the corresponding {\em lower} trend formula
akin to (\ref{lowertrend}).
For example, with $k=1$ and $q=2$, the most probable values of maximal prime gaps $G(x)$
turn out to be less than the EVT-predicted value $\frac{x\log\li x}{\li x}$---less by approximately $\log x \log\log x$ 
(cf.~Cadwell~\cite[p.\,912]{cadwell}).
In this respect, primes do not behave like ``random darts''.
Instead, the situation looks as if primes ``conspire together'' so that each prime $p_n\le x$
lowers the typical maximal gap $G(x)$ by about $p_n^{-1} \log x$; indeed, we have
$\sum_{p_n\le x} p_n^{-1} \sim \log\log x.$
Below we offer a heuristic explanation of this~phenomenon.

Let $\tau_{q,r}(d,x)$ be the number of gaps of a given even size $d=p'-p$ between successive primes
$p,p'\equiv r$ (mod $q$), $p'\le x$.
Empirically, the function $\tau_{q,r}$ has the form (cf.\,\cite{AresCastro2006,goldstonledoan2011,wolf2011})
\begin{equation}\label{defTau}
\tau_{q,r}(d,x) ~\approx~ P_q(d) B_q(x) e^{-d\,\cdot A_q(x)},
\end{equation}

\noindent
where $P_q(d)$ is an oscillating factor (encoding a form of {\em singular series}), and
\begin{equation}\label{tauZero}
\tau_{q,r}(d,x)~=~P_q(d)~=~0 \quad\mbox{ if } q\nmid d \mbox{ \  or \ } 2\nmid d.
\end{equation}

The essential point now is that
we can find the unknown functions $A_q(x)$ and $B_q(x)$ in (\ref{defTau})
just by assuming the exponential decay of $\tau_{q,r}$ as a function of $d$
and employing the following two conditions (which are true by definition of
$\tau_{q,r}$):
\begin{equation}\label{condpi}
\mbox{(a)~~the total number of gaps is }
\sum_{d=2}^{G_{q,r}(x)} \tau_{q,r}(d,x) ~\approx~ \pi(x;q,r);  \qquad\qquad\qquad
\end{equation}
\begin{equation}\label{condx}
\mbox{(b)~~the total length of gaps is~~\, }
\sum_{d=2}^{G_{q,r}(x)} d \cdot \tau_{q,r}(d,x) ~\approx~ x.   \qquad\qquad\qquad\qquad
\end{equation}

The erratic behavior of the oscillating factor $P_q(d)$ presents an obstacle in the calculation of sums~(\ref{condpi}) and (\ref{condx}).
We will assume that, for sufficiently regular functions $f(d,x)$,
\begin{equation}\label{avg1}
\sum_{d} P_q(d) f(d,x) ~\approx~ s \sum_{d}f(d,x),
\end{equation}

\noindent
where $s$ is such that, on average, $P_q(d) \approx s$;
and the summation is for $d$ such that both sides of (\ref{avg1}) are non-zero.
Extending the summation in Equations~(\ref{condpi}), (\ref{condx}) to infinity, using (\ref{avg1}),
and writing
%
$$d=cj, \quad j\in{\mathbb N}, \quad c=\lcm(2,q)=O(q), \quad s = \lim\limits_{n\to\infty}\frac{1}{n}\sum_{j=1}^n P_q(cj),
$$

\noindent
we obtain two series expressions:
(\ref{condpi}) gives us a geometric series
\begin{equation}\label{seriespi}
\sum_{d=2}^{\infty} \tau_{q,r}(d,x)
~\approx~ {sB_q(x)} \sum_{j=1}^\infty e^{-cjA_q(x)}
~=~ {sB_q(x)} \cdot \frac{e^{-cA_q(x)} }{1-e^{-cA_q(x)}}
~\approx~ \pi(x;q,r),
\end{equation}

\noindent
while (\ref{condx}) yields a differentiated geometric series
\begin{equation}\label{seriesx}
\sum_{d=2}^{\infty} d \cdot \tau_{q,r}(d,x)
~\approx~ {csB_q(x)} \sum_{j=1}^\infty je^{-cjA_q(x)}
~=~ {cs B_q(x)} \cdot \frac{e^{-cA_q(x)}}{(1-e^{-cA_q(x)})^2}
~\approx~ x.
\end{equation}

Thus we have obtained two equations:
$$
{sB_q(x)} \cdot \frac{e^{-cA_q(x)}}{1-e^{-cA_q(x)}} ~\approx~ \pi(x;q,r),
\qquad
{csB_q(x)}\cdot \frac{e^{-cA_q(x)}}{(1-e^{-cA_q(x)})^2} ~\approx~ x.
$$

To solve these equations, we use the approximations $e^{-cA_q(x)}\approx1$
and $1-e^{-cA_q(x)}\approx cA_q(x)$ (which is justified because we expect $A_q(x)\to0$ for large $x$).
In this way we obtain
\begin{equation}\label{solutionAB}
A_q(x) \approx \frac{\pi(x;q,r)}{x}, \qquad
B_q(x) \approx \frac{c \pi^2(x;q,r)}{sx}.
\end{equation}
A posteriori we indeed see that $A_q(x)\to0$ as $x\to\infty$.
Substituting (\ref{solutionAB}) into (\ref{defTau}) we get
\begin{equation}\label{solutionTau}
\tau_{q,r}(d,x) ~\approx~ P_q(d) \, \frac{c \pi^2(x;q,r)}{sx} e^{-d\cdot\pi(x;q,r)/x}.
\end{equation}

From (\ref{solutionTau}) we can obtain an approximate formula for $G_{q,r}(x)$.
Note that $\tau_{q,r}(d,x)=1$ when the gap of size $d$ is maximal---in which case we have $d = G_{q,r}(x)$.
So, to get an approximate value of the maximal gap $G_{q,r}(x)$,
we solve for $d$ the equation $\tau_{q,r}(d,x)=1$, or
\begin{equation}\label{taueq1}
\frac{c \pi^2(x;q,r)}{x} e^{-d\cdot\pi(x;q,r)/x} ~\approx~ 1,
\end{equation}

\noindent
where we skipped $P_q(d)/s$ because, on average, $P_q(d)\approx s$.
Taking the log of both sides of (\ref{taueq1}) we find the solution $G_{q,r}(x)$
expressed directly in terms of $\pi(x;q,r)$:
\begin{equation}\label{solutionGqr}
G_{q,r}(x) ~\approx~ \frac{x}{\pi(x;q,r)} \cdot
\left( \log \frac{\pi^2(x;q,r)}{x} + \log c \right).
\end{equation}

Since $\pi(x;q,r)\approx \displaystyle\frac{\li x }{\varphi(q)}$
and $\displaystyle\log\frac{\pi^2(x;q,r)}{x} \approx 2\log\frac{\li x}{\varphi(q)} - \log x$,
we can state the following

\vspace{6pt}
\bfprop{Conjecture on the trend of $G_{q,r}(x)$.}
The most probable sizes of maximal gaps $G_{q,r}(x)$ are near a trend curve $T(q,x)$:
\begin{equation}\label{GqrTrend}
G_{q,r}(x) ~\sim~ T(q,x) ~=~
\frac{\varphi(q) x}{\li x} \cdot \left( 2\log\frac{\li x}{\varphi(q)} - \log x + b \right),
\end{equation}

\noindent
where $b=b(q,x)=O(\log q)$ tends to a constant as $x\to\infty$.
The difference $G_{q,r}(x)- T(q,x)$ changes its sign infinitely often.

Further, we expect that the {\em width of distribution} of the maximal gaps near $x$ is $O_q(\log x)$;
i.e., the~width of distribution is on the order of the average gap $\varphi(q) \log x$ (see Section~\ref{numDistrib}).
On the other hand, for large $x$, the trend (\ref{GqrTrend}) differs from the line $\varphi(q) \log^2 x$
by $O_q(\log x \log\log x)$, that is, by~much more than the average gap.
This suggests natural generalizations of the Cram\'er and Shanks conjectures:

\bfprop{Generalized Cram\'er conjecture for $G_{q,r}(p)$.}
Almost all maximal gaps $G_{q,r}(p)$ satisfy
\begin{equation}\label{GqrCramer}
G_{q,r}(p) ~<~ \varphi(q) \log^2 p.
\end{equation}

\bfprop{Generalized Shanks conjecture for $G_{q,r}(p)$.}
Almost all maximal gaps $G_{q,r}(p)$ satisfy
\begin{equation}\label{GqrShanks}
G_{q,r}(p) ~\sim~ \varphi(q) \log^2 p \qquad\mbox{ as }p\to\infty.
\end{equation}

\smallskip\noindent
Conjectures (\ref{GqrCramer}) and (\ref{GqrShanks})
can be viewed as particular cases of (\ref{GcCramer}), (\ref{GcShanks}) for $k=1$.

\subsection{How Many Maximal Gaps Are There?}\label{howmanymaxgaps}
This section generalizes the heuristic reasoning of \cite[Sect.\,2.3]{kourbatov2017}.
Let $R_c(n)$ be the size of the $n$-th record (maximal) gap between primes in ${\mathbb P}_c$.
Denote by $N_c(x)$ the total number of maximal gaps observed between primes  in ${\mathbb P}_c$ not exceeding $x$.
Let $\ell=\ell(x;q,\cH)$ be a continuous slowly varying function estimating $\mean\limits_r(N_c(ex)-N_c(x))$,
the average number of maximal gaps between primes in ${\mathbb P}_c$,
with the upper endpoints $p'\in[x,ex]$.
For $x\to\infty$, we will heuristically argue that if the limit of $\ell$ exists, then the limit is $k+1$.
Suppose that
$$
\displaystyle\lim_{x\to\infty} \frac{N_c(x)}{\log x}
~=~\lim_{x\to\infty}\frac{\mean_r N_c(x)}{\log x}
~=~\lim_{x\to\infty}\ell(x;q,\cH)~=~\ell_* > 0,
$$

\noindent
and the limit $\ell_*$ is independent of $q$.
Let $n$ be a ``typical'' number of maximal gaps up to $x$;
our~assumption $\lim\limits_{x\to\infty}\ell = \ell_*$ means that
\begin{equation}\label{nvslogx}
n ~\sim~ \ell_* \log x \qquad\mbox{ as }x\to\infty.
\end{equation}

For large $n$, we can estimate the order of magnitude of the typical $n$-th maximal gap $R_c(n)$
using the generalized Cram\'er and Shanks conjectures (\ref{GcCramer}) and (\ref{GcShanks}):
\begin{equation}\label{rough}
R_c(n)~=~G_{c}(x)
~\lesssim~ C_{k,\cH}^{-1}\, \varphi_{k,\cH}(q) \log^{k+1} x
~\sim~  C_{k,\cH}^{-1}\, \varphi_{k,\cH}(q) \frac{n^{k+1}}{\ell_*^{\,k+1}}.
\end{equation}

Define $\Delta R_c(n) = R_c(n+1)-R_c(n).$
By formula (\ref{rough}), for large $q$ and large $n$ we have
$$
\mean_r R_c(n)
~\sim~ C_{k,\cH}^{-1}\, \varphi_{k,\cH}(q) \frac{n^{k+1}}{\ell_*^{\,k+1}},                            
$$
\begin{align*}
\mean_r \Delta R_c(n)
&~=~ \mean_r \big( R_c(n+1)-R_c(n) \big)                                               \\
&~\sim~ \frac{C_{k,\cH}^{-1}\, \varphi_{k,\cH}(q)}{\ell_*^{\,k+1}} \cdot \big((n+1)^{k+1}-n^{k+1}\big) \\
&~\sim~ \frac{C_{k,\cH}^{-1}\, \varphi_{k,\cH}(q)}{\ell_*^{\,k+1}} \cdot (k+1)n^k,
\end{align*}

\noindent
where the mean is taken over all $\cH$-allowed residue classes; see Sect.\,\ref{defphik}.
Combining this with (\ref{nvslogx}), we find
\begin{equation}\label{meanRsimKp1}
\mean_r \Delta R_c(n) ~\sim~ \frac{k+1}{\ell_*} \cdot C_{k,\cH}^{-1}\, \varphi_{k,\cH}(q)\log^{k} x.
\end{equation}

On the other hand, heuristically we expect that, on average, two consecutive record gaps
should differ by the ``local'' average gap (\ref{bara}) between primes in ${\mathbb P}_c$:
\begin{equation}\label{avgGap}
\mean_r \Delta R_c(n) ~\sim~ C_{k,\cH}^{-1}\, \varphi_{k,\cH}(q)\log^{k} x \ \ \mbox{ ($~\sim~$average gap near $x$)}.
\end{equation}

Together, Equations (\ref{meanRsimKp1}) and (\ref{avgGap}) imply that
$$
\ell_* ~=~ k+1.
$$

Therefore, for large $x$ we should expect (see Sections~\ref{countingMaxGaps} and \ref{irtime}; cf.\,\cite{krug2007})
\begin{equation}\label{Ncvslogx}
N_c(x) ~\sim~ (k+1)\log x \qquad\mbox{ as }x\to\infty.
\end{equation}

\noindent
{\bf Special cases. } For the number $N_{q,r}$ of maximal gaps between primes $p\equiv r$ (mod $q$) we have
\begin{equation}\label{Nqrvslogx}
N_{q,r}(x) ~\sim~ 2\log x \qquad\mbox{ as }x\to\infty.
\end{equation}
This is asymptotically equivalent to the following semi-empirical formula for the number of maximal prime gaps up to $x$
(i.e., for the special case $k=1$, $q=2$; see \cite[Sect.\,3.4; OEIS \seqnum{A005669}]{kourbatov2016}):
\begin{equation}\label{twologlix}
N_{2,1}(x) ~\sim~ 2\log \li x \qquad\mbox{ as }x\to\infty.
\end{equation}

Formula (\ref{twologlix}) tells us that maximal prime gaps occur, on average,
about {\em twice} as often as records in an i.i.d.~random sequence of $\lfloor\li x\rfloor$ terms.
Note also the following straightforward generalization of~({\ref{twologlix}) giving
a very rough estimate of $N_{q,r}(x)$ for $k=1$:
\begin{equation}\label{twologlixphiq}
N_{q,r}(x) ~\approx~ \max\left(0,\,2\log\frac{\li x }{\varphi(q)}\right).
\end{equation}

Computation shows that, for the special case of maximal prime gaps $G(x)$,
formula (\ref{twologlix}) works quite well.
However, the more general formula ({\ref{twologlixphiq}) usually {\it overestimates} $N_{q,r}(x)$.
At the same time, the~right-hand side of ({\ref{twologlixphiq}) is {\it less than} $2\log x$.
Thus the right-hand sides of (\ref{Nqrvslogx}) as well as (\ref{twologlixphiq})
overestimate the actual gap counts $N_{q,r}(x)$ in most cases.

In Section \ref{countingMaxGaps} we will see an alternative (a posteriori) approximation
based on the average number of maximal gaps observed for primes in the interval $[x,ex]$.
Namely, the estimated average number $\ell(x;q,\cH)$ of maximal gaps with endpoints in $[x,ex]$ is
%
\begin{equation}\label{ellHyperbola}
\ell(x;q,\cH)~\approx~ \mean_r (N_c(ex)-N_c(x))
~\approx~ k+1 - \frac{\kappa(q,\cH)}{\log x + \delta(q,\cH)}.
\end{equation}

\section{Numerical Results}\label{numresults}
To test our conjectures of the previous section, we performed extensive computational experiments.
We used PARI/GP (see Appendix \ref{AppendixA} for code examples) 
to compute maximal gaps $G_c$ between initial primes $p=r+nq\in{\mathbb P}_c$
in densest admissible prime $k$-tuples, $k\le6$.
We experimented with many different values of $q\in[4,10^5]$.
To assemble a complete data set of maximal gaps for a given $q$,
we~used all $\cH$-allowed residue classes $r$ (mod $q$).
For additional details of our computational experiments with maximal gaps between primes $p=r+nq$
(i.e., for the case $k=1$), see also \cite[Sect.\,3]{kourbatov2016}.
In~this section we omit the subscript $\cH$ in $\varphi_k$ and $C_k$ because
we are working with {\it densest} $k$-tuples:
for each $k=2,4,6$ there is only one densest pattern $\cH$, while
for each $k=3,5,7$ there are two densest patterns~$\cH$,
with equal numerical values of functions $\varphi_k(q)$ and equal Hardy--Littlewood constants $C_k$.

\pagebreak

\subsection{The Growth Trend of Maximal Gaps}\label{numTrend}

The vast majority of maximal gap sizes $G_c(x)$ are indeed observed near the trend curves predicted
in Section \ref{maxgaps}.
Specifically, for maximal gaps $G_c$ between primes $p=r+nq\in{\mathbb P}_c$ in $k$-tuples ($k\ge2$),
the gap sizes are mostly found in the neighborhood 
of the corresponding trend curves of Eqs.\,(\ref{lowertrend}), (\ref{uppertrend}) derived from extreme value theory.
However, for $k=1$, the trend Eq.\,(\ref{GqrTrend}) gives a better prediction
of maximal gaps $G_{q,r}$. 

Figures \ref{figqr313trend}--\ref{figse313trend} illustrate our numerical results for
$k=1,2,6$, $q=313$.
The horizontal axis in these figures is $\log^{k+1}{\negthinspace}p$ for {\em end-of-gap} primes $p$.
Note that all gaps shown in the figures satisfy the generalized Cram\'er conjecture, i.e.,
inequalities (\ref{GcCramer}), (\ref{GqrCramer}); for rare exceptions, see Section~\ref{ExtraLargeGaps}.
Results for other values of $q$ look similar to Figures~\ref{figqr313trend}--\ref{figse313trend}.

\vspace*{1.6cm}

\begin{figure}[H] 
\centering
\includegraphics[bb=2 3 579 356,width=5.67in,height=3.9in,keepaspectratio]{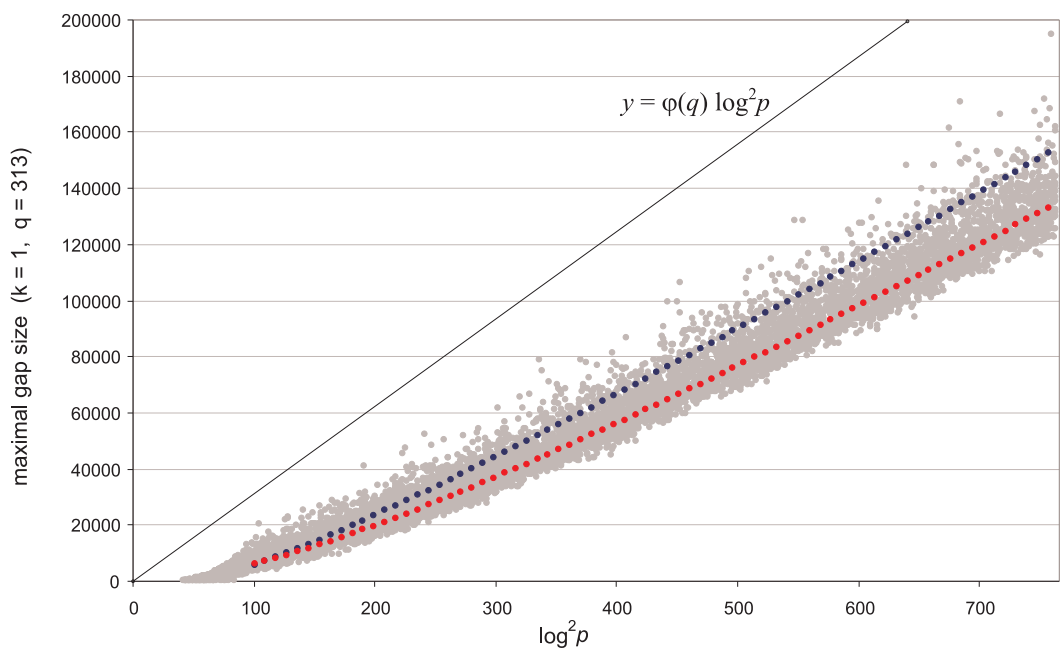}
\caption{Maximal gaps $G_{q,r}$ between primes $p=r+nq\le x$ for $q=313$, $x<10^{12}$. 
Red curve: trend~(\ref{GqrTrend}),\,(\ref{defb});
\ blue curve: EVT-based trend $\frac{\varphi(q) x }{\li x} \log\frac{\li x}{\varphi(q)}$;
\ top line: $y = \varphi(q) \log^2 p$.
}
\label{figqr313trend}
\end{figure}
\unskip
\begin{figure}[H] 
\centering
\includegraphics[bb=2 3 579 356,width=5.67in,height=3.8in,keepaspectratio]{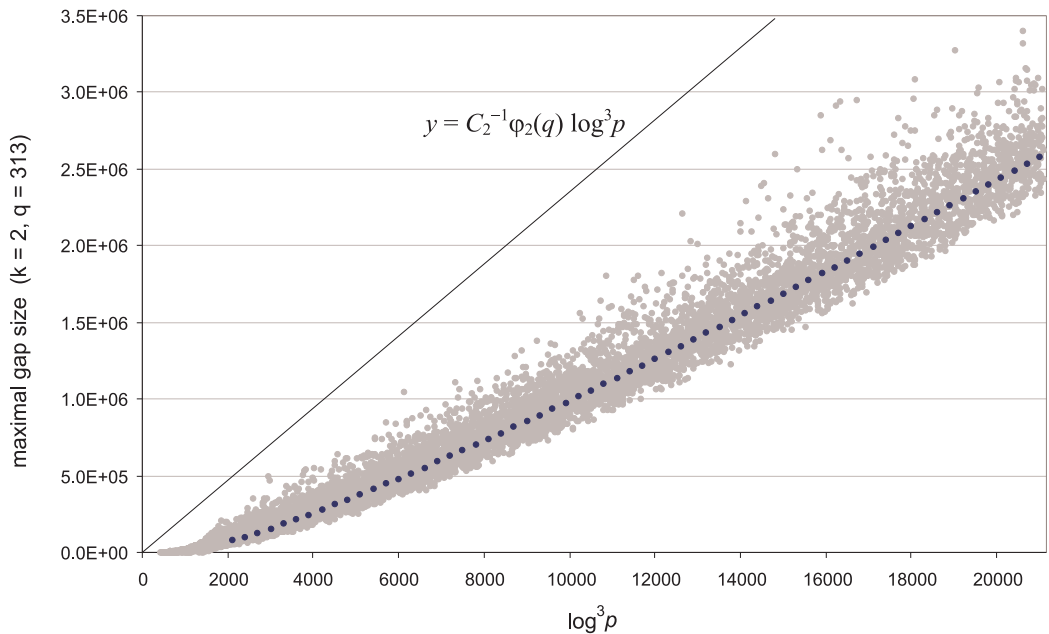}
\caption{Maximal gaps $G_{c}$ between lesser twin primes $p=r+nq\in{\mathbb P}_c$ below $x$
for $q=313$, $x<10^{12}$, $k=2$.
Dotted curve: trend $ T_c$ of Equation~(\ref{lowertrend});
top line: $y = C_{2}^{-1}\varphi_{2}(q)\log^3 p$.
}
\label{figtw313trend}
\end{figure}

Numerical evidence suggests that

\begin{itemize}
\item
For $k=1$ (the case of maximal gaps $G_{q,r}$ between primes $p=r+nq$) the EVT-based trend curve
$\frac{\varphi(q) x }{\li x} \log\frac{\li x}{\varphi(q)}$
goes too high (Fig.\,\ref{figqr313trend}, blue curve). Meanwhile, the trend (\ref{GqrTrend}) 
$$
T(q,x) ~=~
\frac{\varphi(q) x}{\li x} \cdot \left( 2\log\frac{\li x}{\varphi(q)} - \log x + b \right) \quad\mbox{(Fig.\,\ref{figqr313trend}, red curve)}
$$
satisfactorily predicts gap sizes $G_{q,r}(x)$, with the empirical correction term

\begin{equation}\label{defb}
b ~=~ b(q,x) ~\approx\, \left(b_0+\frac{b_1}{(\log\log x)^{b_2}}\right)\log\varphi(q) ~\asymp~ \log\varphi(q),
\end{equation}

\noindent
where the parameter values

\begin{equation}\label{b012}
b_0=1, \qquad b_1=4, \qquad b_2=2.7
\end{equation}

\noindent
are close to optimal for $q\in[10^2,10^5]$ and $x\in[10^7,10^{14}]$.
Here the qualifier {\em optimal} is to be understood in conjunction with
the rescaling transformation (\ref{defw}) introduced below in Section~\ref{numDistrib}.
A trend $T(q,x)$ is optimal if after transformation (\ref{defw}) the most probable rescaled values
$w$ turn out to be near zero, and the mode of best-fit Gumbel distribution for $w$-values
is also close to zero, $\mu\approx0$; see Figure~\ref{figqr16001hist}.
In view of (\ref{defb}) it is possible that, for all $q$,
the optimal term $b$ in (\ref{GqrTrend}) has the form
$b(q,x) = (1 + \beta(q,x))\cdot\log\varphi(q) \sim \log\varphi(q)$,
where $\beta(q,x)$ very slowly decreases to zero as $x\to\infty$.
(Note that in Section~\ref{caseOfPrimes} we correctly estimated $b$ to be $O(\log q)$
but did not predict the appearance of Euler's function $\varphi(q)$ in the term $b$.)

\item
For $k=2$, approximately half of maximal gaps $G_c$ between lesser twin primes $p\in{\mathbb P}_c$
are below the {\em lower} trend curve $ T_c(x)$ of Equation~(\ref{lowertrend}),
while the other half are above that curve; see Figure~\ref{figtw313trend}.

\item
For $k\ge3$, more than half of maximal gaps $G_c$ are usually above
the lower trend curve $ T_c(x)$ of Equation~(\ref{lowertrend}).
At the same time, more than half of maximal gaps are usually below
the upper trend curve $\bar T_c(x)$ of Equation~(\ref{uppertrend}); see Figure~\ref{figse313trend}.
Recall that the two trend curves $ T_c$ and $\bar T_c$ are within $k\bar a_c$ from each other
as $x\to\infty$; see (\ref{trenddiff}).
\end{itemize}
\vspace{-18pt}

\begin{figure}[H] 
\centering

\includegraphics[bb=2 3 579 356,width=5.67in,height=3.8in,keepaspectratio]{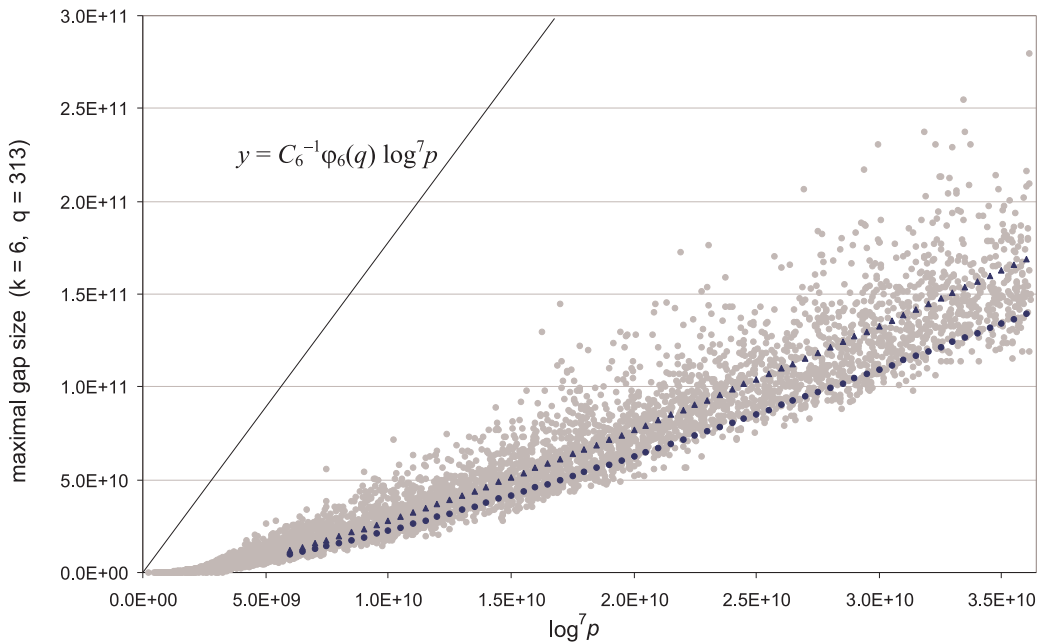}
\caption{Maximal gaps $G_{c}$ between prime sextuplets $p=r+nq\in{\mathbb P}_c$ below $x$
for $q=313$, $x<10^{14}$, $k=6$.
Dotted curves: trends $ T_c$ ($\bullet$) and $\bar T_c$ ($\blacktriangle$) of Equations~(\ref{lowertrend}) and (\ref{uppertrend});
top line: $y = C_{6}^{-1}\varphi_{6}(q)\log^7 p$.
}
\label{figse313trend}
\end{figure}

\medskip
As noted by Brent \cite{brent2014}, twin primes seem to be more random than primes.
We can add that, likewise, maximal gaps $G_{q,r}$ between primes in a residue class
seem to be somewhat less random than those for prime $k$-tuples;
primes $p\equiv r$ (mod $q$) do not go quite as far from each other
as we would expect based on extreme value theory.
Pintz \cite{pintz2007} discusses various other aspects of the ``random'' and not-so-random behavior of primes.

\subsection{The Distribution of Maximal Gaps}\label{numDistrib}

In Section \ref{numTrend} we have tested equations that determine the growth trend
of maximal gaps between primes in sequences ${\mathbb P}_c$.
How are maximal gap sizes distributed in the neighborhood of their respective~trend?

We will perform a rescaling transformation (motivated by extreme value theory):
subtract the trend from the actual gap size, and then divide the result by a natural unit,
the ``local'' average gap. This way each maximal gap size is mapped to its rescaled value:
$$
\mbox{ maximal gap size } G ~~\mapsto~
\mbox{ rescaled value} = \frac{G ~-~ \mbox{trend}}{\mbox{average gap}}.
$$

Gaps above the trend curve are mapped to positive rescaled values, while
gaps below the trend curve are mapped to negative rescaled values.



\medskip\noindent
{\bf Case }$k=1$. For maximal gaps $G_{q,r}$ between primes $p\equiv r$ (mod $q$),
the trend function $T$ is given by Equations~(\ref{GqrTrend}), (\ref{defb}) and (\ref{b012}).
The rescaling operation has the form
\begin{equation}\label{defw}
G_{q,r}(x) ~\mapsto~
w = \frac{G_{q,r}(x) - T(q,x)}{a(q,x)}.
\end{equation}

\noindent
where $a(q,x)=\displaystyle\frac{\varphi(q)x}{\li x}$.
Figure \ref{figqr16001hist} shows histograms of rescaled values $w$ for maximal gaps $G_{q,r}$
between primes $p \equiv r$ (mod $q$) for $q=16001$.

\medskip\noindent
{\bf Case }$k\ge2$. For maximal gaps $G_c$ between prime $k$-tuples with $p=r+nq\in{\mathbb P}_c$,
we can use the trend $ T_c$ of Equation~(\ref{lowertrend}).
Then the rescaling operation has the form
\begin{equation}\label{defh}
G_{c}(x) ~\mapsto~
h = \frac{G_{c}(x) -  T_c(x)}{ a_c(x)},
\end{equation}

\noindent
where $ a_c(x)$ is defined by (\ref{tildea}).
Figure \ref{figtw16001hist} shows histograms of rescaled values $ h$ for maximal gaps $G_{c}$
between lesser twin primes $p=r+nq\in{\mathbb P}_c$ for $q$ = 16001, $k=2$.

\begin{figure}[H] 
\centering
\includegraphics[bb=7 4 351 356,width=5.2in,height=5.2in,keepaspectratio]{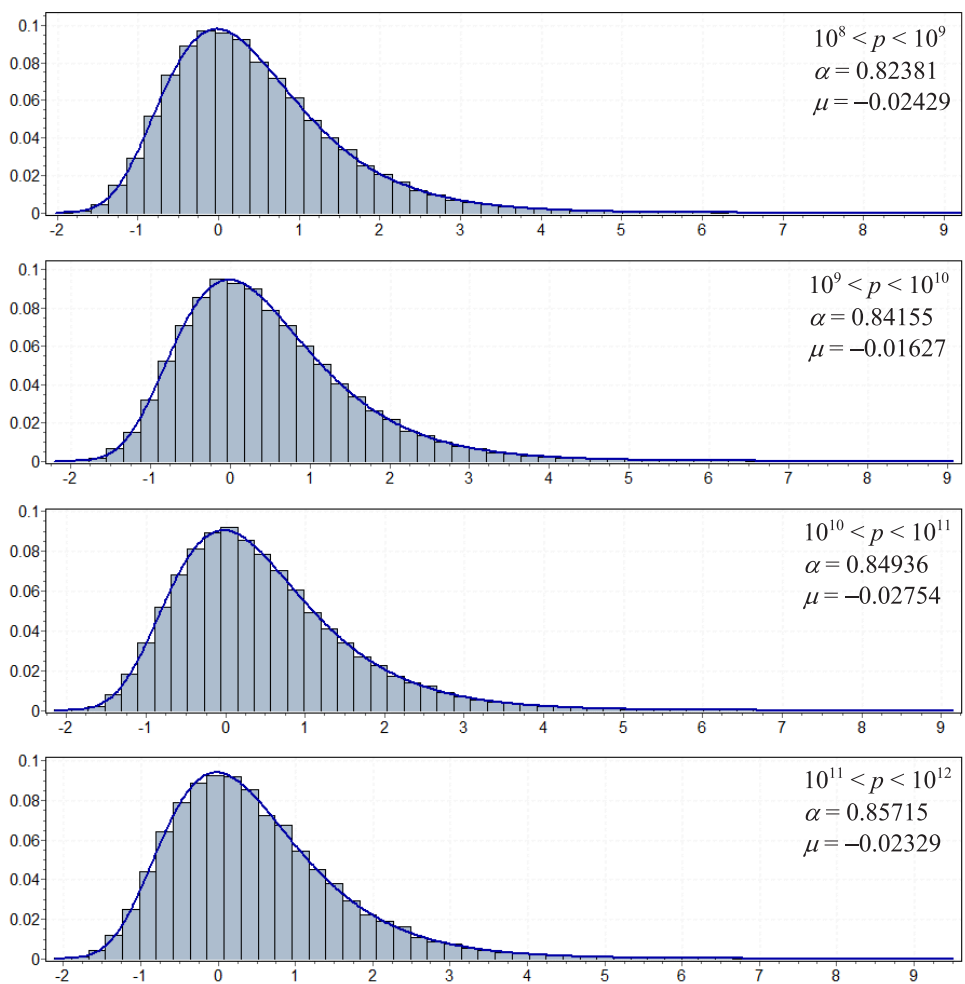}
\caption{Histograms of $w$-values (\ref{defw}) for maximal gaps $G_{q,r}$ between primes
$p=r+nq$ for $q=16001$, $r\in[1,16000]$. Curves are best-fit Gumbel distributions (pdfs)
with scale $\alpha$ and mode $\mu$.}
\label{figqr16001hist}
\end{figure}
In both Figures~\ref{figqr16001hist} and \ref{figtw16001hist},
note that the histograms and fitting distributions are skewed to the right,
i.e., the right tail is longer and heavier.
Among two-parameter distributions, the Gumbel extreme value distribution
is a very good fit; cf.\,\cite{kourbatov2014,liprattshakan}.
This was true in all our computational experiments.

\begin{figure}[H] 
\centering
\includegraphics[bb=7 4 351 356,width=5.2in,height=5.2in,keepaspectratio]{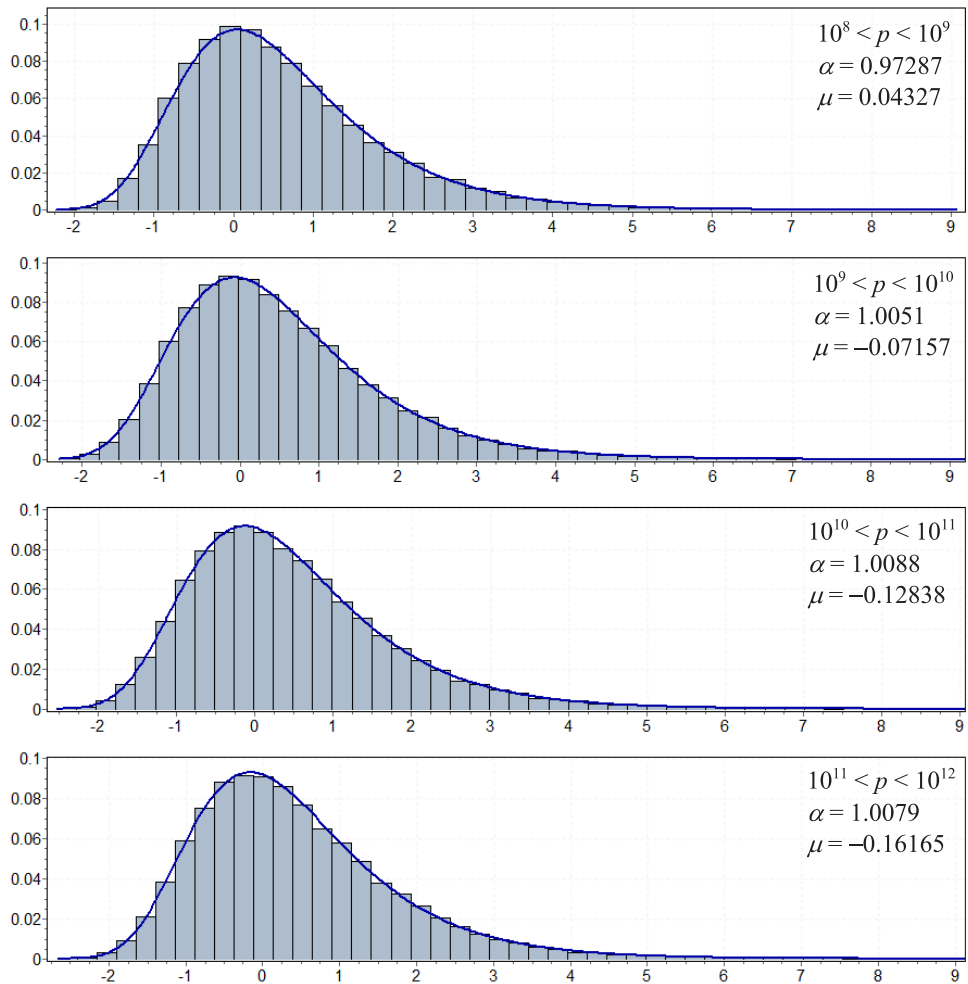}
\caption{Histograms of $ h$-values (\ref{defh}) for maximal gaps $G_{c}$ between lesser twin primes
$p=r+nq\in{\mathbb P}_c$ for $q=16001$ and $\cH$-allowed residue classes $r\in[1,16000]$, $r\ne15999$.
Curves are best-fit Gumbel distributions (pdfs) with scale $\alpha$ and mode $\mu$.
}
\label{figtw16001hist}
\end{figure}

\begin{Remark}For all histograms shown in Figures~\ref{figqr16001hist} and \ref{figtw16001hist},
the Kolmogorov--Smirnov goodness-of-fit statistic is less than 0.01;
in fact, for most of the histograms, the goodness-of-fit statistic is about 0.003.
\end{Remark}

If we look at three-parameter distributions, then an excellent fit is the
{\em Generalized Extreme Value} (GEV) distribution, which includes the Gumbel distribution as a special case.
The {\em shape parameter} in the best-fit GEV distributions is close to zero;
note that the Gumbel distribution is a GEV distribution whose shape parameter is exactly zero.
So could the Gumbel distribution be the limit law for appropriately rescaled sequences of maximal gaps
$G_{q,r}(p)$ and $G_{c}(p)$ as $p\to\infty$? Does such a limiting distribution exist at all?

\pagebreak

\noindent{\bf The scale parameter }$\alpha$.
For $k=1$, we observed that the scale parameter of best-fit Gumbel distributions for $w$-values (\ref{defw})
was in the range $\alpha\in[0.7,1]$. The parameter $\alpha$ seems to slowly grow towards~1 as $p\to\infty$;
see Figure~\ref{figqr16001hist}.
For $k\ge2$, the scale parameter of best-fit Gumbel distributions for $ h$-values~(\ref{defh})
was usually a little over 1; see Figure~\ref{figtw16001hist}.
However, if instead of (\ref{defh}) we use the (simpler) rescaling~transformation
\begin{equation}\label{defbarh}
G_{c}(x) ~\mapsto~
\bar h = \frac{G_{c}(x) - \bar T_c(x)}{\bar a_c(x)},
\end{equation}

\noindent
where $\bar a_c$ and $\bar T_c$ are defined, respectively, by (\ref{bara}) and (\ref{uppertrend}),
then the resulting Gumbel distributions of $\bar h$-values will typically have scales $\alpha$ a little below 1.
In a similar experiment with {\em random} gaps, the~scale was also close to 1; see
\cite[Sect.\,3.3]{kourbatov2016}.

\subsection{Counting the Maximal Gaps}\label{countingMaxGaps}
We used PARI/GP function {\tt findallgaps} (see source code in {Appendix} \ref{PARIauxfunc})
to determine average numbers of maximal gaps $G_{q,r}$ between primes $p=r+nq$, $p\in[x,ex]$,
for $x=e^j$, $j=1,2,\ldots,27$. Similar statistics were also gathered for gaps $G_c$.
Figures \ref{figmeanrecK1}--\ref{figmeanrecK4} show the results of this computation
for $q=16001$, $k\le4$. When $x$ is large, the average number
of maximal gaps $G_c$ for $p\in[x,ex]$ indeed seems to {\em very slowly} approach $k+1$,
as predicted by Equation~(\ref{Ncvslogx}).
When $x$ is small ($x<q/e$), there is {\em at most} one prime $p\in[2,ex]$ in sequence $\Pc$---and often there are no such primes at all;
accordingly, we~see no gaps ending in $[x,ex]$, and the corresponding plot points
in Figures~\ref{figmeanrecK1}--\ref{figmeanrecK4} are zero.

Starting from some $x_0>q/e$, the gap counts in $[x,ex]$ are no longer zero.
Here we observe a ``transition region'' in which the mean number of maximal gaps $G_c(p)$
for primes $p \in [x,ex]$ grows from 0 to a little over 1,
while $x$ increases by about 3 orders of magnitude from $x_0$.
The non-monotonic behavior of plot points in the transition region is explained, in part, by the fact that
here the gap size may be comparable to the size of intervals $[x,ex]$.
Then, for larger $x$, the typical number of gaps $G_c$ between $k$-tuples continues to slowly increase;
specifically, the graph of $\mean(N_c(ex)-N_c(x))$ vs.~$\log x$  is closely approximated
by a hyperbola with horizontal asymptote $y=k+1$; see (\ref{Ncvslogx}), (\ref{ellHyperbola}) in Section \ref{howmanymaxgaps}.

Why do the observed curves resemble hyperbolas?
If we were working with  {\em random} gaps,
then perhaps the curves could be explained using the theory of records; cf.\,\cite{krug2007}.
But primes are not random numbers; and so we simply treat the hyperbolas
in Figures \ref{figmeanrecK1}--\ref{figmeanrecK4} as an {\em experimental~fact}.

\begin{figure}[H] 
\centering\vspace{-18pt}

\includegraphics[bb=2 3 579 356,width=4.9in,height=3.6in,keepaspectratio]{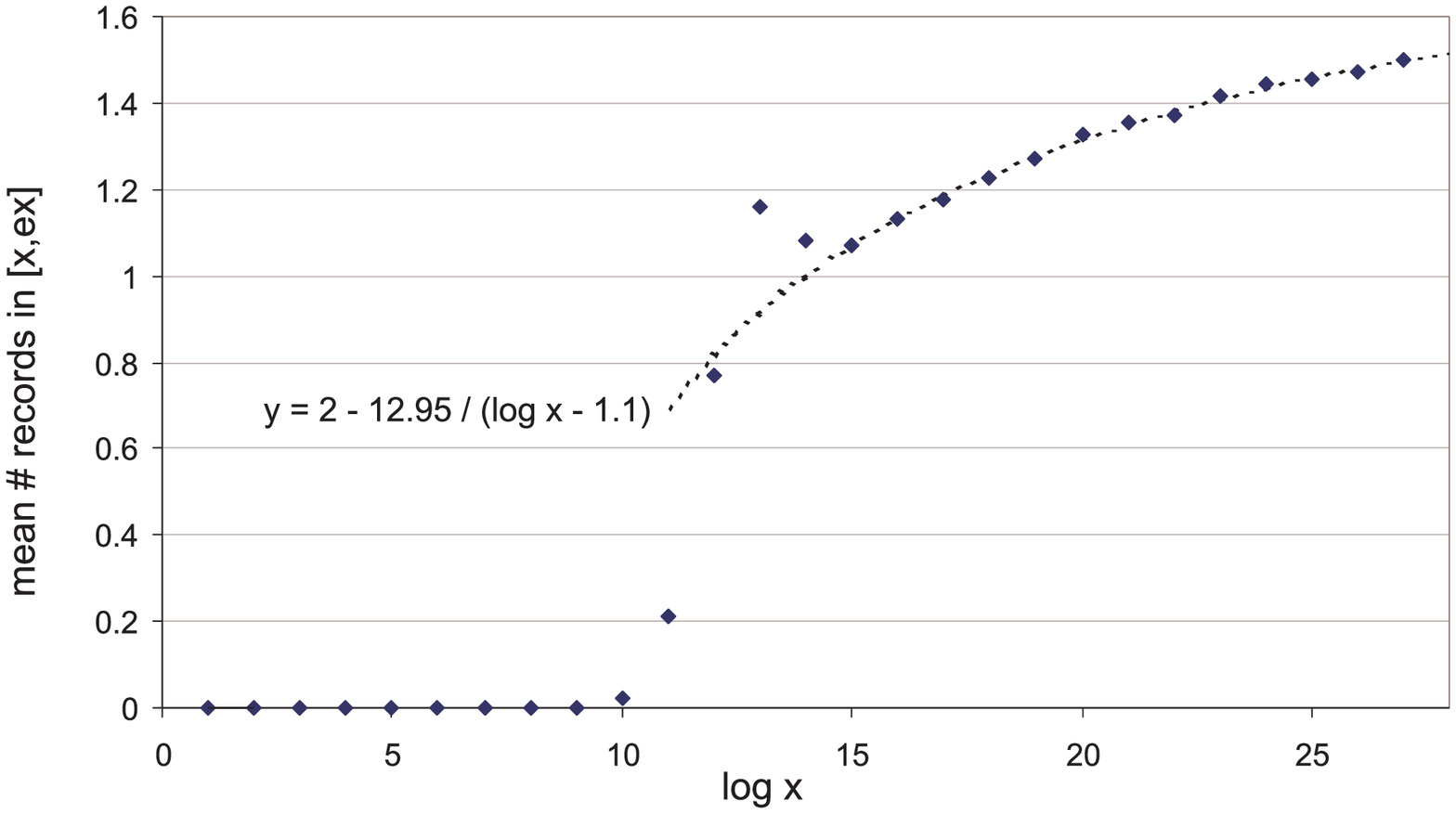}
\caption{Primes $p=r+nq$, $k=1$, $q=16001$.
Mean number of maximal gaps $G_{q,r}$ observed for $p\in[x,ex]$, $x=e^j$, $j\le27$.
Averaging for all $\cH$-allowed $r$. Dotted curve is a hyperbola with horizontal asymptote $y=2$.
}
\label{figmeanrecK1}

\vspace*{1.6cm}

\includegraphics[bb=2 3 579 356,width=4.9in,height=3.6in,keepaspectratio]{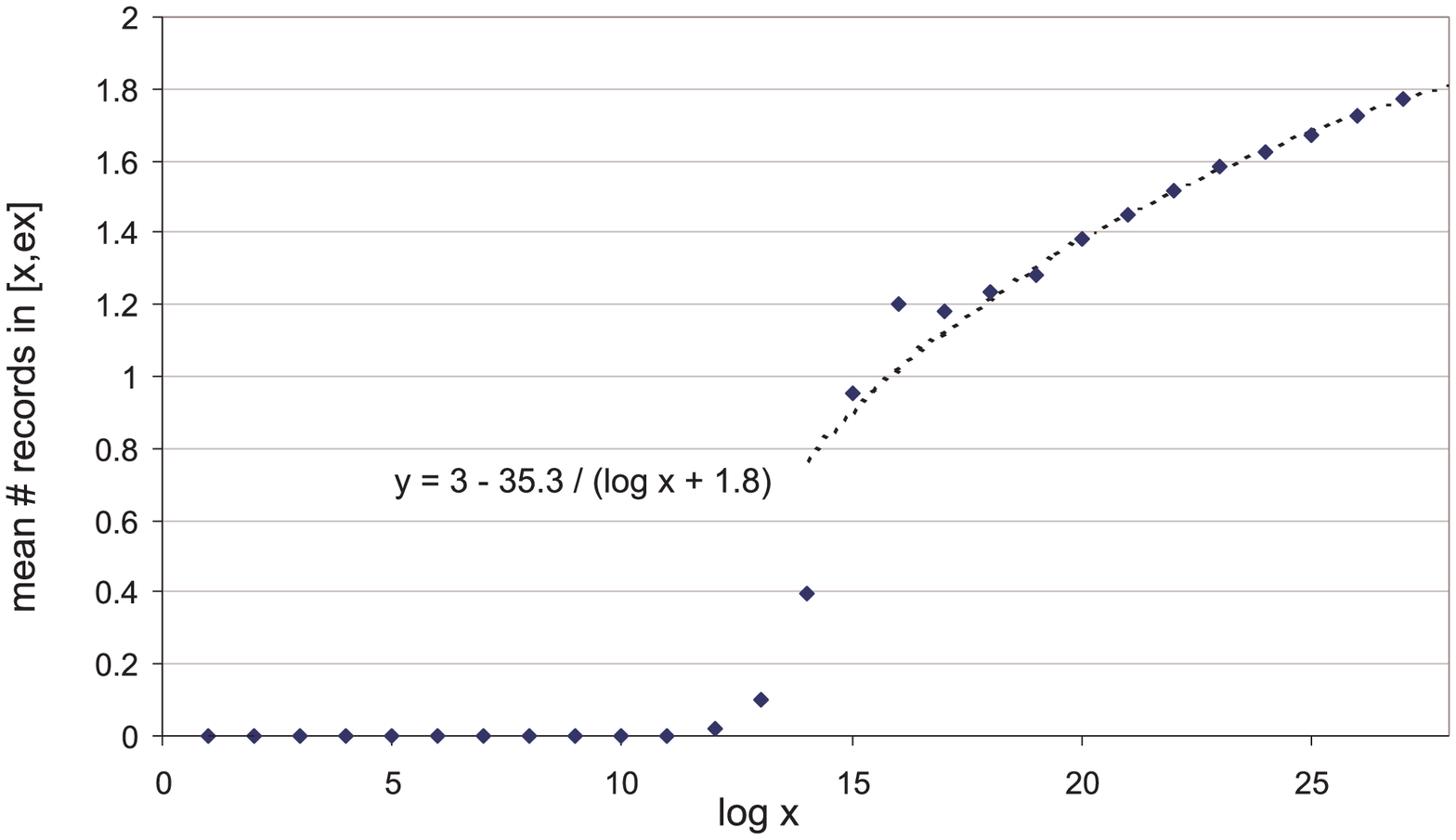}
\caption{Lesser twin primes $p=r+nq\in{\mathbb P}_c$, $k=2$, $q=16001$.
Mean number of maximal gaps $G_{c}$ observed for $p\in[x,ex]$, $x=e^j$, $j\le27$.
Averaging for all $\cH$-allowed $r$. Dotted curve is a hyperbola with horizontal asymptote $y=3$.
}
\label{figmeanrecK2}
\end{figure}

\begin{figure}[H] 
\centering
\includegraphics[bb=2 3 579 356,width=5.0in,height=3.6in,keepaspectratio]{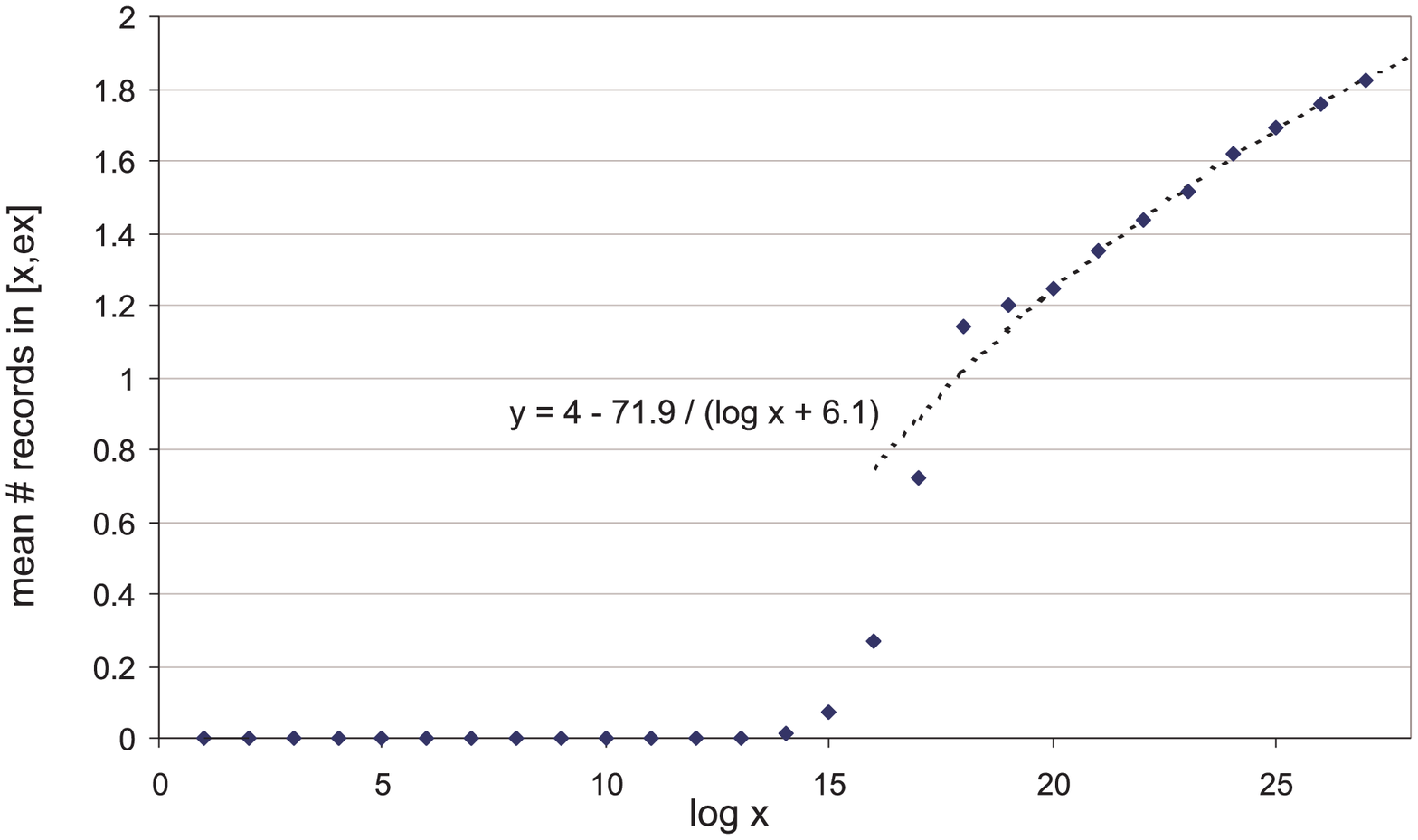}
\caption{Prime triplets ($p$, $p+2$, $p+6$), $p=r+nq\in{\mathbb P}_c$, $k=3$, $q=16001$.
Mean number of maximal gaps $G_{c}$ observed for $p\in[x,ex]$, $x=e^j$, $j\le27$.
Averaging for all $\cH$-allowed $r$. Dotted curve is a hyperbola with horizontal asymptote $y=4$.
}
\label{figmeanrecK3}

\vspace*{1.6cm}

\includegraphics[bb=2 3 579 356,width=4.9in,height=3.6in,keepaspectratio]{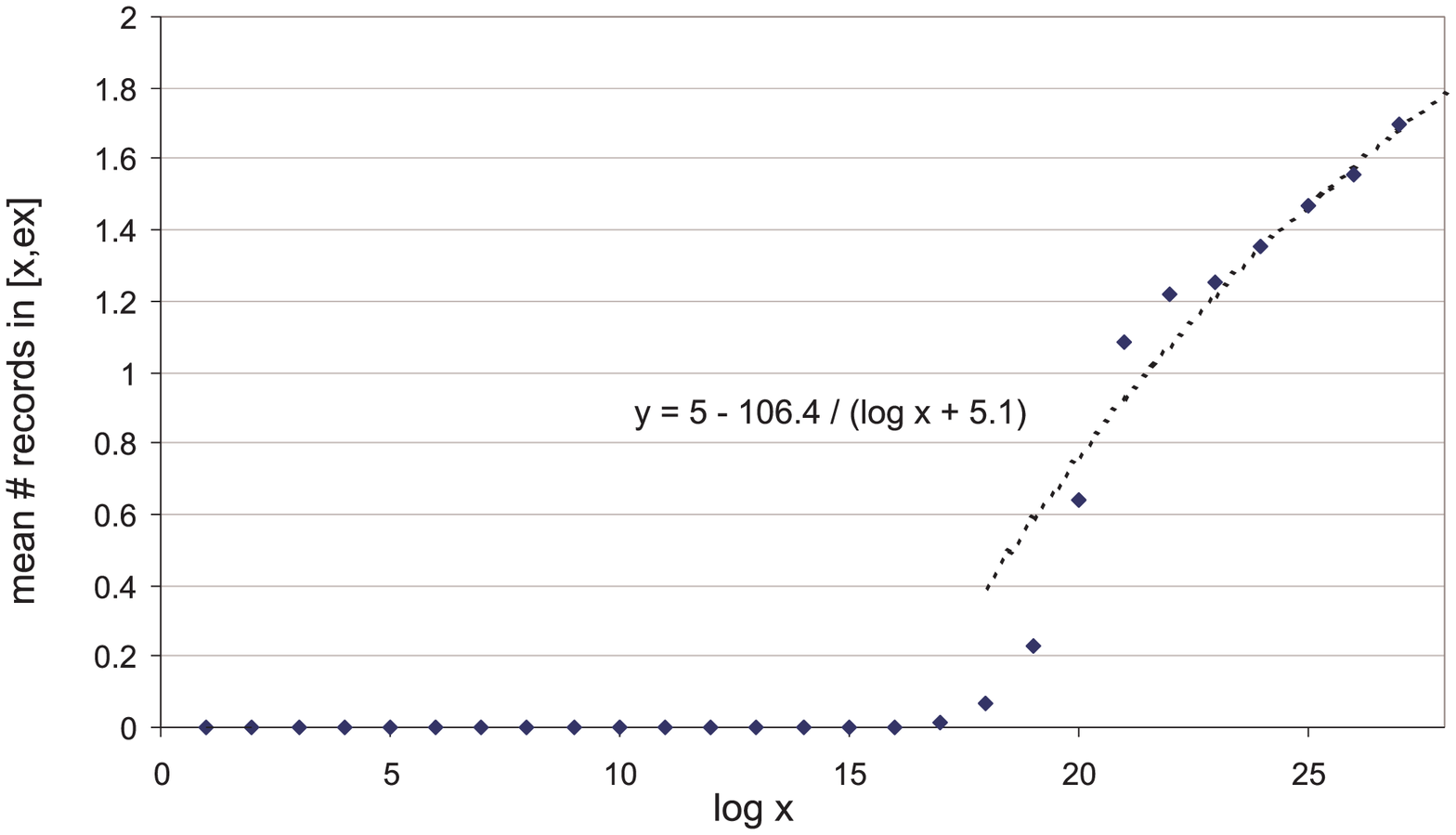}
\caption{Prime quadruplets ($p$, $p+2$, $p+6$, $p+8$), $p=r+nq\in{\mathbb P}_c$, $k=4$, $q=16001$.
Mean number of maximal gaps $G_{c}$ observed for $p\in[x,ex]$, $x=e^j$, $j\le27$.
Averaging for all $\cH$-allowed $r$.
Dotted curve is a hyperbola with horizontal asymptote $y=5$.
}
\label{figmeanrecK4}
\end{figure}

\pagebreak
\subsection{How Long Do We Wait for the Next Maximal Gap?}\label{irtime}

Let $P(n)=\mbox{\seqnum{A002386}}(n)$ and $P'(n)=\mbox{\seqnum{A000101}}(n)$
be the lower and upper endpoints of the $n$-th record (maximal) gap $R(n)$ between primes:
$R(n) = \mbox{\seqnum{A005250}}(n) = P'(n)-P(n)$.

Consider the distances $P(n)-P(n-1)$ from one maximal gap to the next.
(In statistics, a similar quantity is sometimes called ``inter-record times'').
In Figure~\ref{figirtimes} we present a plot of these distances;
the~figure also shows the corresponding plot for {\em twin primes}.
As can be seen from Figure~\ref{figirtimes}, the~quantity $P(n)-P(n-1)$ grows
approximately exponentially with $n$ (but not monotonically).
Indeed, typical inter-record times are expected to satisfy\footnote{
 The asymptotic equivalence $\sim$ in Eqs.\,(\ref{irtimeprimes}) and (\ref{irtimektuples})
 is a restatement of Eqs.\,(\ref{Ncvslogx})~and~(\ref{Nqrvslogx}).
 It would be logically unsound to suppose that $\log(P(n)-P(n-1)) \stackrel{?}{\sim} \log P(n)$
 because we cannot exclude the possibility that $\log(P(n)-P(n-1))$
 might (very rarely) become as small as \mbox{$\log G(x)\approx 2\log\log x$}, where $x=P(n)$.
}
\begin{equation}\label{irtimeprimes}
\log(P(n)-P(n-1)) ~<~ \log P(n) ~\sim~ \frac{n}{2} \qquad\mbox{ as }n\to\infty.
\end{equation}

\begin{figure}[H] 
\centering
\includegraphics[width=5.6in,height=4.3in,keepaspectratio]{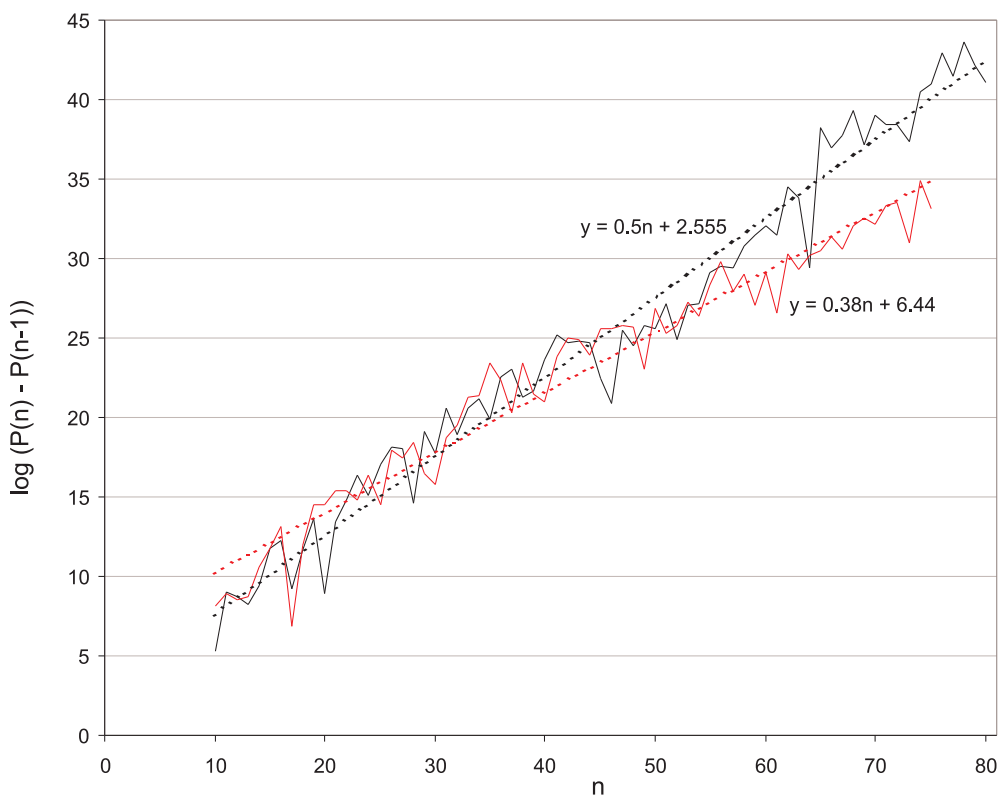}
\caption{Inter-record times $P(n)-P(n-1)$ for gaps between primes (black) and
a similar quantity $P_c(n)-P_c(n-1)$ for gaps between twin primes (red).
Lines are exponential fits. Values for $n<10$ are~skipped.
}
\label{figirtimes}
\end{figure}

More generally, let $P_c(n)$ and $P'_c(n)$ be the endpoints of the $n$-th maximal gap $R_c(n)$
between primes in sequence ${\mathbb P}_c$, where each prime is $r$ (mod $q$)
and starts an admissible prime $k$-tuple.
Then, in accordance with heuristic reasoning of Section~\ref{howmanymaxgaps},
for typical inter-record times $P_c(n)-P_c(n-1)$ separating the maximal gaps $R_c(n-1)$ and $R_c(n)$
we expect to see
\begin{equation}\label{irtimektuples}
\log(P_c(n)-P_c(n-1)) ~<~ \log P_c(n) ~\sim~ \frac{n}{k+1} \qquad\mbox{ as }n\to\infty.
\end{equation}

In the special case $k=2$, that is, for maximal gaps between {\em twin primes},
the right-hand side of (\ref{irtimektuples}) is expected to be
$\frac{n}{3}$ for large $n$
(whereas Figure~\ref{figirtimes} suggests the right-hand side $0.38n$
based on a very limited data set for $10\le n\le75$).
As we have seen in Section~\ref{countingMaxGaps},
the average number of maximal gaps between $k$-tuples occurring for primes $p\in[x,ex]$
slowly approaches $k+1$ {\em from below}.
For moderate values of $x$ attainable in computation, this average is typically between 1 and $k+1$.
Accordingly, we see that the right-hand side of (\ref{irtimektuples})
yields a prediction $\asymp e^{n/(k+1)}$ that {\em underestimates} the typical inter-record times and the primes $P_c(n)$.
Computations may yield estimates
$$P_c(n)-P_c(n-1) < P_c(n) \approx Ce^{\beta n},
$$
where
$$\beta\in[\frac{1}{k+1},1],
$$

\noindent
with the estimated value of $\beta$ depending on the range of available data.

\begin{Remark} Sample graphs of $\log P_c(n)$ vs.~$n$ can be plotted online at the OEIS website:
click {\em graph} and scroll to the logarithmic plot for sequences
\seqnum{A002386} ($k=1$),
\seqnum{A113275} ($k=2$),
\seqnum{A201597} ($k=3$),
\seqnum{A201599} ($k=3$),
\seqnum{A229907} ($k=4$),
\seqnum{A201063} ($k=5$),
\seqnum{A201074} ($k=5$),
\seqnum{A200504} ($k=6$).
In all these graphs, when $n$ is large enough,
$\log P_c(n)$ seems to grow approximately linearly with $n$.
We conjecture that the slope of such a linear approximation slowly decreases,
approaching the slope value $1/(k+1)$ as $n\to\infty$.
\end{Remark}

\subsection{Exceptionally Large Gaps: $G_{q,r}(p)>\varphi(q)\log^2p$}\label{ExtraLargeGaps}

Recall that for the {\it maximal prime gaps} $G(x)$
Shanks \cite{shanks} conjectured the asymptotic equality
$G(x)\sim \log^2 x$, a strengthened form of Cram\'er's conjecture.
This seems to suggest that (unusually large) maximal gaps $g$
may in fact occur as early as at $x \asymp e^{\sqrt{g}}$.
On the other hand, Wolf \cite{wolf1997} conjectured that typically a gap of size $d$
appears for the first time between primes near $\sqrt{d}\cdot e^{\sqrt{d}}$.
Combining these observations, we may further observe that exceptionally large maximal gaps, that is,
\begin{equation}\label{XLgaps}
\mbox{ largest gaps } g ~=~ G(x) ~>~ \log^2 x
\end{equation}
are also those which appear for the first time unusually early.
Namely, they occur at $x$ roughly by a factor of $\sqrt d$ earlier than
the typical first occurrence of a gap $d$ at $x\asymp \sqrt{d}\cdot e^{\sqrt{d}}$.
Note that Granville~\cite[p.\,24]{granville} suggests that gaps of unusually large size (\ref{XLgaps})
occur infinitely often---and we will even see infinitely many of those exceeding $1.1229 \log^2 x$.
In contrast, Sun \cite[Conj.\,2.3]{sun2013}
made a conjecture implying that exceptions like (\ref{XLgaps}) occur only finitely often,
while Firoozbakht's conjecture implies that exceptions (\ref{XLgaps}) never occur for primes
$p\ge 11$; see \cite{kourbatov2015u}.
Here we cautiously predict that exceptional gaps of size (\ref{XLgaps})
are only a zero proportion of maximal gaps. This can be viewed as restatement
of the generalized Cram\'er conjectures (\ref{GcCramer}), (\ref{GqrCramer})
for the special case $k=1$, $q=2$.

Table \ref{XLgapsTable} lists exceptionally large maximal gaps $G_{q,r}(p)$ between primes $p\equiv r$ (mod $q$)
for which inequality (\ref{GqrCramer}) does not hold:

$$\mbox{ largest gaps } G_{q,r}(p)>\varphi(q)\log^2p.$$

\pagebreak


\begin{table}[H]
\caption{Exceptionally large maximal gaps: $G_{q,r}(p)>\varphi(q)\log^2p$ \ for $p<10^9$, \ $r<q\le30000$.  \label{XLgapsTable}}
\centering
\begin{tabular}{rrrrrc}
\toprule { \large$\vphantom{1^{1^1}}$}
\textbf{Gap \boldmath{$G_{q,r}(p)$}}  &  \textbf{Start of Gap}  &  \textbf{End of Gap (\boldmath{$p$})} &   \boldmath{$q~~$} & \boldmath{$r~~$} & \boldmath{$G_{q,r}(p)/(\varphi(q)\log^2p$)} \\
[0.5ex]\midrule
\vphantom{\fbox{$1^1$}}
(i)~~~~~208650
&   3415781 & 3624431 & 1605 & 341 & 1.0786589153 \\
316790  &    726611 & 1043401 & 2005 & 801 & 1.0309808771 \\
229350  &   1409633 & 1638983 & 2085 & 173 & 1.0145547849 \\
532602  &    355339 &  887941 & 4227 & 271 & 1.0081862161 \\
984170  &   5357381 & 6341551 & 4279 &  73 & 1.0339720553 \\
1263426 &  10176791 & 11440217 & 4897 & 825 & 1.0056800570 \\
2306938 &  82541821 & 84848759 & 6907 & 3171 & 1.0022590147 \\
3415794 & 376981823 & 380397617 & 8497 & 3921 & 1.0703375544 \\
2266530 & 198565889 & 200832419 & 8785 & 7319 & 1.0335372951 \\
7326222 & 222677837 & 230004059 & 20017 & 8729 & 1.0166221904 \\
6336090 &  10862323 & 17198413 & 23467 & 20569 & 1.0064940453 \\
7230930 & 130172279 & 137403209 & 24595 & 15539 & 1.0468373915 \\
5910084 & 51763573 & 57673657 & 28971 & 21367 & 1.0199911211 \\  \midrule
(ii)~~~~411480
& 470669167 & 471080647 & 3048 & 55 & 1.0235488825 \\
208650  &   3415781 & 3624431 & 3210 & 341 & 1.0786589153 \\
316790  &    726611 & 1043401 & 4010 & 801 & 1.0309808771 \\
229350  &   1409633 & 1638983 & 4170 & 173 & 1.0145547849 \\
657504  & 896016139 & 896673643 & 4566 & 2563 & 1.0179389550 \\
1530912 & 728869417 & 730400329 & 6896 & 3593 & 1.0684247390 \\
532602  &    355339 & 887941 & 8454 & 271 & 1.0081862161 \\
984170  &   5357381 & 6341551 & 8558 & 73 & 1.0339720553 \\
1263426 &  10176791 & 11440217 & 9794 & 825 & 1.0056800570 \\
2119706 & 665152001 & 667271707 & 10046 & 6341 & 1.0223668231 \\
1885228 & 163504573 & 165389801 & 10532 & 5805 & 1.0000704209 \\
1594416 & 145465687 & 147060103 & 13512 & 9007 & 1.0026889378 \\
2306938 &  82541821 & 84848759 & 13814 & 3171 & 1.0022590147 \\
3108778 & 524646211 & 527754989 & 15622 & 12585 & 1.0098218219 \\
1896608 &    164663 & 2061271 & 16934 & 12257 & 1.0598397341 \\
3415794 & 376981823 & 380397617 & 16994 & 3921 & 1.0703375544 \\
2266530 & 198565889 & 200832419 & 17570 & 7319 & 1.0335372951 \\
2937868 &  71725099 & 74662967 & 17698 & 12803 & 1.0103309882 \\
2823288 &  37906669 & 40729957 & 18098 & 9457 & 1.0162761199 \\
2453760 &  11626561 & 14080321 & 18176 & 12097 & 1.0107626289 \\
3906628 & 190071823 & 193978451 & 18692 & 11567 & 1.1480589845 \\
2157480 & 13074917 & 15232397 & 27660 & 19397 & 1.0716522452 \\
5450496 & 366870073 & 372320569 & 28388 & 11949 & 1.0140771094 \\
3422630 & 735473 & 4158103 & 29762 & 21185 & 1.0368176014 \\  \midrule
(iii)~~~657504
& 896016139 & 896673643 & 2283 & 280 & 1.0179389550 \\
2119706 & 665152001 & 667271707 & 5023 & 1318 & 1.0223668231 \\
3108778 & 524646211 & 527754989 & 7811 & 4774 & 1.0098218219 \\
1896608 &    164663 & 2061271 & 8467 & 3790 & 1.0598397341 \\
2937868 &  71725099 & 74662967 & 8849 & 3954 & 1.0103309882 \\
2823288 &  37906669 & 40729957 & 9049 & 408 & 1.0162761199 \\
3422630 &    735473 & 4158103 & 14881 & 6304 & 1.0368176014 \\
3758772 & 144803717 & 148562489 & 15927 & 11360 & 1.0000152764 \\
3002682 &   8462609 & 11465291 & 16869 & 11240 & 1.0107025944  \\
8083028 & 344107541 & 352190569 & 19619 & 9900 & 1.1134625422  \\
4575906 &  20250677 & 24826583 & 22653 & 21548 & 1.0463153374 \\
5609136 & 34016537 & 39625673 & 26967 & 11150 & 1.0412524005 \\
7044864 & 302145839 & 309190703 & 27519 & 14738 & 1.0048671503 \\
6580070 & 9659921 & 16239991 & 28609 & 18688 & 1.0046426332 \\
\bottomrule
\end{tabular}
\end{table}

\pagebreak

Three sections of Table \ref{XLgapsTable} correspond to (i) odd $q,r$; (ii) even $q$; (iii) even $r$.
(Overlap between sections is due to the fact that $\varphi(q)=\varphi(2q)$ for odd $q$.)
No other maximal gaps with this property were found for $p<10^9$, $q\le30000$.
No such large gaps exist for $p<10^{10}$, $q\le1000$.

\begin{Remark}It is interesting that, for every gap listed in Table \ref{XLgapsTable},
at least one of the numbers $q$ and $r$ is {\it composite}.
Thus far we have never seen a gap violating (\ref{GqrCramer}) with both $q$ and $r$ prime.
\end{Remark}

\section{Summary}

We have extensively studied record (maximal) gaps between prime $k$-tuples in residue classes (mod $q$).
Our computational experiments described in Section \ref{numresults} took months of computer time.
Numerical evidence allows us to arrive at the following conclusions, which are also supported by
heuristic reasoning.

\begin{itemize}
\item For $k=1$, the observed growth trend of maximal gaps $G_{q,r}(x)$ is given by (\ref{GqrTrend}) and (\ref{defb}).
In~particular, for maximal prime gaps ($k=1$, $q=2$) the trend equation reduces to
$$
G_{2,1}(x) ~{\sim}~T(2,x) = \log^2 x - 2\log x\log\log x + O(\log x).
$$

\item For $k\ge2$, a significant proportion of maximal gaps $G_{c}(x)$
are observed between the trend curves of Equations~(\ref{lowertrend}) and (\ref{uppertrend}),
which can be heuristically derived from extreme value theory.

\smallskip
\item The Gumbel distribution, after proper rescaling, is a possible limit law for $G_{q,r}(p)$ as well as $G_{c}(p)$.
The existence of such a limiting distribution is an open question.

\smallskip
\item Almost all maximal gaps $G_{q,r}(p)$ between primes in residue classes mod $q$
seem to satisfy appropriate generalizations of the Cram\'er and Shanks conjectures
(\ref{GqrCramer}) and (\ref{GqrShanks}):
$$G_{q,r}(p) ~\lesssim~ \varphi(q) \log^{2} p.$$

\smallskip
\item Similar generalizations (\ref{GcCramer}) and (\ref{GcShanks}) of the Cram\'er and Shanks conjectures
are apparently true for almost all maximal gaps $G_{c}(p)$ between primes in ${\mathbb P}_c$:
$$G_{c}(p) ~\lesssim~ C_{k,\cH}^{-1}\, \varphi_{k,\cH}(q) \log^{k+1} p.$$

\smallskip
\item Exceptionally large gaps $G_{q,r}(p)>\varphi(q)\log^2p$ are extremely rare (Table \ref{XLgapsTable}).
We conjecture that only a zero proportion of maximal gaps are such exceptions.
A similar observation holds for $G_c(p)$ violating (\ref{GcCramer}).

\smallskip
\item We conjecture that the total number $N_{q,r}(x)$ of maximal gaps $G_{q,r}$ observed up to $x$ is
below $C\log x$ for some $C>2$.

\smallskip
\item More generally, we conjecture: the number $N_{c}(x)$ of maximal gaps between primes in ${\mathbb P}_c$ up to $x$
satisfies the inequality $N_{c}(x)<C\log x$ for some $C>k+1$, where $k$ is the number of integers in
the pattern $\cH$ defining the sequence ${\mathbb P}_c$.
\end{itemize}





\acknowledgments{We are grateful to the anonymous referees for useful suggestions. Thanks also to all
contributors and editors of the websites {\it OEIS.org} and {\it PrimePuzzles.net}.}




\pagebreak

\appendixtitles{yes}
\appendix
\section{Details of Computational Experiments}\label{AppendixA}

Interested readers can reproduce and extend our results using the programs below.

\subsection{PARI/GP Program maxgap.gp}

\begin{verbatim}
default(realprecision,11)
outpath = "c:\\wgap"

\\ maxgap(q,r,end [,b0,b1,b2]) ver 2.1 computes maximal gaps g
\\ between primes p = qn + r, as well as rescaled values (w, u, h):
\\   w - as in eqs.(33),(45)-(47) of arXiv:1901.03785 (this paper); 
\\   u - same as w, but with constant b = ln phi(q);
\\   h - based on extreme value theory (cf. randomgap.gp in arXiv:1610.03340)
\\ Results are written on screen and in the folder specified by outpath string.
\\ Computation ends when primes exceed the end parameter.
maxgap(q,r,end,b0=1,b1=4,b2=2.7) = {
  re = 0;
  p = pmin(q,r);
  t = eulerphi(q); 
  inc = q;
  while(p<end,
    m = p + re;
    p = m + inc;
    while(!isprime(p), p+=inc);
    while(!isprime(m), m-=inc);
    g = p - m;
    if(g>re,
      re=g; Lip=li(p); a=t*p/Lip; Logp=log(p);
      h = g/a-log(Lip/t);
      u = g/a-2*log(Lip/t)+Logp-log(t);
      w = g/a-2*log(Lip/t)+Logp-log(t)*(b0+b1/max(2,log(Logp))^b2);
      f = ceil(Logp/log(10));
      write(outpath"\\"q"_1e"f".txt",
            w" "u" "h" "g" "m" "p" q="q" r="r);
      print(w" "u" "h" "g" "m" "p" q="q" r="r);
      if(g/t>log(p)^2, write(outpath"\\"q"_1e"f".txt","extra large"));
      if(g%2==0, inc=lcm(2,q));

      \\ optional part: statistics for p in intervals [x/e,ex] for x=e^j
      i = ceil(Logp);
      j = floor(Logp);
      if(N!='N,N[j]++);  \\ count maxima with p in [x,ex] for x=e^j
      write(outpath"\\"q"_exp"i".txt", w" "u" "h" "g" "m" "p" q="q" r="r);
      write(outpath"\\"q"_exp"j".txt", w" "u" "h" "g" "m" "p" q="q" r="r);
    )
  )
}
\end{verbatim}


\subsection{PARI/GP: Auxiliary Functions for maxgap.gp}\label{PARIauxfunc}

\begin{verbatim}
\\ These functions are intended for use with the program maxgap.gp
\\ It is best to include them in the same file with maxgap.gp

\\ li(x) computes the logarithmic integral of x
li(x) = real(-eint1(-log(x)))

\\ pmin(q,r) computes the least prime p = qn + r, for n=0,1,2,3,...
pmin(q,r) = forstep(p=r,1e99,q, if(isprime(p), return(p)))

\\ findallgaps(q,end): Given q, call maxgap(q,r,end) for all r coprime to q.
\\ Output total and average counts of maximal gaps in intervals [x,ex].
findallgaps(q,end) = {
  t = eulerphi(q);
  N = vector(99,j,0);
  for(r=1,q, if(gcd(q,r)==1,maxgap(q,r,end)));
  nmax = floor(log(end));
  for (n=1,nmax,
    avg = 1.0*N[n]/t;
    write(outpath"\\"q"stats.txt", n" "avg" "N[n]);
  )
}
\end{verbatim}

\subsection{Notes on Distribution Fitting}

In order to study distributions of rescaled maximal gaps,
we used the distribution-fitting software EasyFit \cite{easyfit}.
Data files created with {\tt maxgap.gp} 
are easily imported into EasyFit:

\begin{enumerate}

\item From the {\it File} menu, choose {\it Open}.

\item Select the data file.

\item Specify {\it Field Delimiter} = {\it space}.

\item Click {\it Update}, then {\it OK}.
\end{enumerate}

\medskip\noindent
{\bf Caution:} PARI/GP outputs large and small real numbers in a mantissa-exponent format
{\em with a space} preceding the exponent (e.g.,\ {\tt1.7874829515 E-5}),
whereas EasyFit expects such numbers {\em without a space} (e.g.,\ {\tt1.7874829515E-5}).
Therefore, before importing into EasyFit, search the data files for {\tt" E"}
and replace all occurrences with {\tt"E"}.

\section{The Hardy--Littlewood Constants \boldmath$C_{k,\cH}$}\label{hlconst}
The Hardy--Littlewood $k$-tuple conjecture \cite{hl1923} allows one to predict the average frequencies
of prime $k$-tuples near $p$, as well as the approximate total counts of prime $k$-tuples below $x$.
Specifically, the Hardy--Littlewood $k$-tuple constants $C_{k,\cH}$, divided by $\log^k p$,
give us an estimate of the average frequency of prime $k$-tuples near $p$:
$$
\mbox{Frequency of $k$-tuples } ~\sim~ \frac{C_{k,\cH}}{\log^k p}.
$$
%
Accordingly \cite[pp.\,61--68]{riesel}, for a given $k$-tuple pattern $\cH$,
the total count of $k$-tuples below $x$ is
$$
\pi_{k,\cH}(x) ~\sim~ C_{k,\cH} \int_2^x{\negthinspace}{dt\over\log^k t} ~=~ C_{k,\cH} \,{\Li}_k(x).
$$

The Hardy--Littlewood constants $C_{k,\cH}$ can be defined in terms of infinite products over primes.
In~particular, for densest admissible prime $k$-tuples with $k\le7$ we have:
\begin{align*}
C_1 &~=~ 1 \ \mbox{ (by convention, in accordance with the prime number theorem); }
\\
C_{2,\cH} &~=~ 2 \prod_{p>2} {{p(p-2)}\over{(p-1)^2}} ~\approx~ 1.32032363169373914785562422
\qquad \mbox{ (\seqnum{A005597}, \seqnum{A114907});}
\\
C_{3,\cH} &~=~ {9\over2} \prod_{p>3} {{p^2(p-3)}\over{(p-1)^3}} ~\approx~ 2.85824859571922043243013466
\qquad \mbox{ (\seqnum{A065418});}
\\
C_{4,\cH} &~=~ {27\over2} \prod_{p>4} {{p^3(p-4)}\over{(p-1)^4}} ~\approx~ 4.15118086323741575716528556
\qquad \mbox{ (\seqnum{A065419});}
\\
C_{5,\cH} &~=~ {15^4\over2^{11}} \prod_{p>5} {{p^4(p-5)}\over{(p-1)^5}} ~\approx~ 10.131794949996079843988427
\qquad \mbox{ (\seqnum{A269843});}
\\
C_{6,\cH} &~=~ {15^5\over2^{13}} \prod_{p>6} {{p^5(p-6)}\over{(p-1)^6}} ~\approx~ 17.2986123115848886061221077
\qquad \mbox{ (\seqnum{A269846});}
\\
C_{7,\cH} &~=~ {35^6\over3\cdot2^{22}} \prod_{p>7} {{p^6(p-7)}\over{(p-1)^7}} ~\approx~ 53.9719483001296523960730291
\qquad \mbox{ (\seqnum{A271742})}.
\end{align*}

Forbes \cite{forbes} gives values of the Hardy--Littlewood constants up to $k=24$,
albeit with fewer significant digits; see also \cite[p.\,86]{finch}.
Starting from $k=8$, we may often encounter more than one numerical value of $C_{k,\cH}$ for a single $k$.
(If there are $m$ different patterns $\cH$ of densest admissible prime $k$-tuples for the same $k$,
then we typically have $\lceil{m\over2}\rceil$ different numerical values of $C_{k,\cH}$,
depending on the actual pattern $\cH$ of the $k$-tuple; see \cite{forbes}.)

\section{Integrals \boldmath${\Li}_k(x)$}\label{integralsLikx}

Let $k\in{\mathbb N}$ and $x>1$, and let
\begin{align*}
F_k(x) &~=~ \int {\negthinspace}{dx\over\log^k x}
\qquad \mbox{(indefinite integral);}
\\
{\Li}_k(x) &~=~ \int_2^x{\negthinspace}{dt\over\log^k t}
\qquad \mbox{(definite integral).}
\end{align*}

Denote by $\li x$ the conventional logarithmic integral (principal value):
$$
\li x  ~=~ \int_0^x{\negthinspace}{dt\over\log t}
~=~ \int_2^x{\negthinspace}{dt\over\log t} + 1.04516\ldots
$$

In PARI/GP, an easy way to compute $\li x$ is as follows: \
{\tt li(x) = real(-eint1(-log(x)))}.

The integrals $F_k(x)$ and ${\Li}_k(x)=F_k(x)-F_k(2)$ can also be expressed in terms of $\li x$.
Integration by parts gives
$$
\int  \frac{dx}{\log x} = \frac{x}{\log x}  + \frac{x}{\log^2 x} + \frac{2x}{\log^3  x} +
\frac{6x}{\log^4 x } + \cdots + \frac{(k-2)!x}{\log^{k-1} x}+ (k-1)!\int  \frac{dx}{\log^k x}.
$$

Therefore,
\begin{align*}
F_2(x)&~=~ {1\over1!} \left( \li x - {x\over\log x} \right) + C,
\\
F_3(x)&~=~ {1\over2!} \left( \li x - {x\over\log^2 x} (\log x + 1) \right) + C,
\\
F_4(x)&~=~ {1\over3!} \left( \li x - {x\over\log^3 x} (\log^2 x + \log x + 2) \right) + C,
\\
F_5(x)&~=~ {1\over4!} \left( \li x - {x\over\log^4 x} (\log^3 x + \log^2 x + 2 \log x + 6) \right) + C,
\\
F_6(x)&~=~ {1\over5!} \left( \li x - {x \over\log^5 x} (\log^4 x + \log^3 x + 2 \log^2 x + 6 \log x + 24) \right) + C,
\end{align*}

\noindent
and, in general,
$$
F_{k+1}(x)~=~ {1\over k!} \bigg( \li x - {x \over\log^{k} x} \sum_{j=1}^{k} (k-j)! \,\log^{j-1} x \bigg) + C.
$$

Using these formulas we can compute ${\Li}_k(x)$ for approximating $\pi_c(x)$
(the prime counting function for sequence ${\mathbb P}_c$)
in accordance with the $k$-tuple equidistribution conjecture (\ref{ktuplequidist}):
$$
\pi_c(x) ~\approx~ {C_{k,\cH} \over \varphi_{k,\cH}(q)} {\Li}_k(x) ~=~ {C_{k,\cH} \over \varphi_{k,\cH}(q)} (F_k(x)-F_k(2)).
$$

The values of $\li x$, and hence $\Li_k(x)$, can be calculated without (numerical) integration.
For~example, one can use the following rapidly converging series for $\li x$, 
with $n!$ in the denominator and $\log^n x$ in the numerator
(see \cite[]{Prudnikov-et-al-I}, formulas 1.6.1.8--9):
$$
\li x ~=~ \gamma + \log \log x + \sum_{n=1}^{\infty} {\log^{n} x \over n \cdot n!}
\quad\mbox{ for } x > 1.
$$


\reftitle{References}

\bigskip\noindent
{\bf Keywords:} {Cram\'er conjecture; Gumbel distribution; prime gap; prime $k$-tuple; residue class; Shanks conjecture; totient}

\medskip\noindent
{\bf MSC:} {11A41, 11N05}

\medskip\noindent
{\bf Errata:} \url{http://www.javascripter.net/math/errata/err190103785.pdf}

\bigskip\noindent
\footnotesize{Copyright \textcopyright~2019 by the authors. Licensee MDPI, Basel, Switzerland. This article is an open access
article distributed under the terms and conditions of the Creative Commons Attribution
(CC BY) license (http://creativecommons.org/licenses/by/4.0/).
}



\end{document}